# ON THE NUMBER OF ANTIPODAL OR STRICTLY ANTIPODAL PAIRS OF POINTS IN FINITE SUBSETS OF $\mathbb{R}^d$, III

E. Makai, Jr.[1], H. Martini[2], M. H. Nguyễn[3], V. Soltan[4], I. Talata,[5]

[1] Research (partially) supported by Hungarian National Foundation for Scientific Research, grant no-s K75016, K81146

## Abstract

We improve our earlier upper bound on the numbers of antipodal pairs of points among $n$ points in $\mathbb{R}^3$, to $2n^2/5 + O(n^c)$, for some $c > 2$. We prove that the minimal number of antipodal pairs among $n$ points in convex position in $\mathbb{R}^d$, affinely spanning $\mathbb{R}^d$, is $n + d(d-1)/2 - 1$. Let $\underline{sa}_d^s(n)$ be the minimum of the number of strictly antipodal pairs of points among any $n$ points in $\mathbb{R}^d$, with affine hull $\mathbb{R}^d$, and in strictly convex position. The value of $\underline{sa}_d^s(n)$ was known for $d \leq 3$ and any $n$. Moreover, $\underline{sa}_d^s(n) = \lceil n/2 \rceil$ was known for $n \geq 2d$ even, and $n \geq 4d + 1$ odd. We show $\underline{sa}_d^s(n) = 2d$ for $2d + 1 \leq n \leq 4d - 1$ odd, we determine $\underline{sa}_d^s(n)$ for $d = 4$ and any $n$, and prove $\underline{sa}_d^s(2d - 1) = 3(d - 1)$. The cases $d \geq 5$ and $d + 2 \leq n \leq 2d - 2$ remain open, but we give a lower and an upper bound on $\underline{sa}_d^s(n)$ for them, which are of the same order of magnitude, namely $\Theta\left((d-k)d\right)$. We present a simple example of a strictly antipodal set in $\mathbb{R}^d$, of cardinality const $\cdot 1.5874...^d$. We give simple proofs of the following statements: if $n$ segments in $\mathbb{R}^3$ are pairwise antipodal, or strictly antipodal, then

[1] A. Rényi Inst. Math., L. Eötvös Research Network, Pf. 127, H-1364 Budapest, Hungary, http://www.renyi.hu/~makai, e-mail: makai.endre@renyi.hu

[2] Faculty of Mathematics, University of Technology, 09107 Chemnitz, Germany, e-mail: martini@mathematik.tu-chemnitz.de

[3] State University of Moldova, Str. Alexei Mateevici 60, Chişinău, Moldova

[4] Department of Mathematics, George Mason University, Fairfax, Virgina, USA, e-mail: vsoltan@gmu.edu

[5] Óbuda University, M. Ybl Faculty of Architecture, Department of Mathematics, H-1146 Thököly út 74, Budapest, Hungary, e-mail: talata.istvan@renyi.hu

$n \leq 4$, or $n \leq 3$, respectively, and these are sharp. We describe also the cases of equality.

# 1 Introduction

**1.** For basic notation and usual abbreviations we refer to [16]. *Let $X = \{x_1, \ldots, x_n\} \subset \mathbb{R}^d$ (where $x_i \neq x_j$ for $i \neq j$), and let $P = \operatorname{conv} X$ be its convex hull. We always suppose that the affine hull of $X$ is $\mathbb{R}^d$, hence $n \geq d + 1$.* The points $x_i, x_j \in X$ are called *antipodal, or strictly antipodal*, if there are different parallel supporting hyperplanes $H', H''$ of $P = \operatorname{conv}\{x_1, \ldots, x_n\}$ such that $x_i \in H'$ and $x_j \in H''$, or $\{x_i\} = X \cap H'$ and $\{x_j\} = X \cap H''$, respectively. These concepts were introduced by V. Klee ([18], also cf. [16, p. 420, in first edition IN SECOND EDITION?]), or B. Grünbaum [15], respectively.

We write $|\cdot|$, $\operatorname{conv}(\cdot)$, $\operatorname{aff}(\cdot)$, $\operatorname{cl}(\cdot)$, $\operatorname{bd}(\cdot)$, $\operatorname{int}(\cdot)$, $\operatorname{rel int}(\cdot)$ for the cardinality, convex hull, affine hull, closure, boundary, interior, and relative interior (in its affine hull) of a set in $\mathbb{R}^d$, respectively. The norm of a vector is written as $\|\cdot\|$. For $x, y \in \mathbb{R}^d$ we write $[x.y]$ for the segment with endpoints $x, y$. The unit sphere of $\mathbb{R}^d$ is written as $S^{d-1}$. For polytopes $\operatorname{vert}(\cdot)$ denotes their vertex sets. The *support cone of a convex polytope $P \subset \mathbb{R}^d$ at its vertex $x$* is $\cup_{\lambda \geq 0} \lambda(P - x)$. (This is the dual of the normal cone of $P$ at $x$.) A *$k$-fold pyramid over the base a polytope $P$* is a $k$-fold iterated pyramid over the base $P$ (cf. [16]).

We say that a finite subset $X$ of $\mathbb{R}^d$, whose affine hull is $\mathbb{R}^d$, *is in convex position, or strictly convex position*, if $X \subset \operatorname{bd} P$, or $X \subset \operatorname{vert} P \iff X = \operatorname{vert} P$, respectively. For $d \geq 2$ and $n \geq d + 1$ integers we write

$$\mathcal{A}_d(n) := \{X \subset \mathbb{R}^d \mid |X| = n, \operatorname{aff} X = \mathbb{R}^d\}, \text{ and}$$
$$\mathcal{C}_d(n) := \{X \subset \mathbb{R}^d \mid |X| = n, \operatorname{aff} X = \mathbb{R}^d,$$
$$X \text{ is in convex position}\}, \text{ and}$$
$$\mathcal{C}_d^s(n) := \{X \subset \mathbb{R}^d \mid |X| = n, \operatorname{aff} X = \mathbb{R}^d,$$
$$X \text{ is in strictly convex position}\}.$$

We denote by $a(X)$, or $sa(X)$, *the number of antipodal, or strictly antipodal, pairs $\{x_i, x_j\}$ of $X$*, respectively. By Grünbaum [15], $sa(X)$ is half the vertex number of the *difference body* $P - P$ of $P$. For some other interpretations of $sa(X)$, cf. the introduction of [20].

We write

$$\underline{a}_d(n) := \min\{a(X) \mid X \in \mathcal{C}_d(n)\}, \quad \underline{a}_d^s(n) := \min\{a(X) \mid X \in \mathcal{C}_d^s(n)\},$$
$$\overline{a}_d(n) := \max\{a(X) \mid X \in \mathcal{C}_d(n)\}, \quad \overline{a}_d^s(n) := \max\{a(X) \mid X \in \mathcal{C}_d^s(n)\},$$
$$\underline{sa}_d(n) := \min\{sa(X) \mid X \in \mathcal{C}_d(n)\}, \quad \underline{sa}_d^s(n) := \min\{sa(X) \mid X \in \mathcal{C}_d^s(n)\},$$
$$\overline{sa}_d(n) := \max\{sa(X) \mid X \in \mathcal{C}_d(n)\}, \quad \overline{sa}_d^s(n) := \max\{sa(X) \mid X \in \mathcal{C}_d^s(n)\}.$$

Evidently, for $X \in \mathcal{A}_d(d+1)$ we have $a(X) = sa(X) = d(d+1)/2$, so when investigating the above quantities we may restrict our attention to the cases $n \geq d+2$.

**2.** For an overview on lower and upper estimates of $a(X)$ and $sa(X)$, under different hypotheses, cf. [20, paragraphs 1 and 2] and [21]; see also the more recent survey [22]. We remark that [22, p. 227] contains a misprint: namely, for $v_d'(2d-1)$, i.e., our $\underline{sa}_d^s(2d-1)$, the correct value is $3(d-1)$, cf. our Theorem 2.15.

We have to mention that [20, p. 459, lines 12-19] are false. The correct version of the display formula in line 13 was given by Nguyễn and V. Soltan [33] and [26] as follows: $\underline{a}_d^s(n) = n - 1 + d(d-1)/2$. The celebrated Danzer-Grünbaum Theorem [9] states that an antipodal set in $\mathbb{R}^d$ has cardinality at most $2^d$, and equality stands if and only if the set consists of all vertices of a parallelepiped (an affine $d$-cube). A strengthening for the 3-dimensional case of the Danzer-Grünbaum Theorem is given by K. Bezdek-Bisztriczky-K. Böröczky [2]. They call a convex polyhedron edge-antipodal, if for each edge the two endpoints are antipodal, and they prove that even in this case the number of vertices is at most 8, with equality only for a parallelepiped. (It is elementary that for the plane this number is 4.) It seems to be unknown whether this extends to higher dimensions, but in $\mathbb{R}^d$ the upper bound $(d/2+1)^d$ is known, see Swanepoel [34]. Grünbaum [15] proved that a strictly antipodal set in $\mathbb{R}^3$ has at most five points, which is sharp. (It is elementary that for the plane this number is 3.) A strengthening of the Grünbaum theorem is given by Bisztriczky-K. Böröczky [4]: they prove that if an antipodal set in $\mathbb{R}^3$ has six points, then three of these lie on a plane, and the remaining three lie on a parallel plane. They conclude that if this set were strictly antipodal, then small perturbations would preserve strict antipodality, but then the parallel plane condition would be violated. A generalization of [4] is proved by Schürmann-Swanepoel [30]: all 3-dimensional antipodal sets are classified.

Nguyễn-V. Soltan [26] showed that, for $X \in \mathcal{C}_2^s(n)$, $a(X) = n + k$ and $sa(X) = n - k$, where $k$ is the number of parallel pairs of sides of $P = \operatorname{conv} X$.

This strengthens Grünbaum [15], Remark, pp. 9-10, which asserted that the maximum of $a(X)$, for $X \subset \mathbb{R}^2$ ($X \in \mathcal{C}_2^s(n)$ tacitly assumed), is $\lfloor 3n/2 \rfloor$, for all $n \geq 3$. For similar results, cf. also the paper V. Soltan [32].

A metric variant of the question of $a(X)$ and $sa(X)$, for given $n$ and $X \in \mathcal{A}_d(n)$, is treated in Pach-Swanepoel [27] and [28]. This is the question of maximal number of *(strict) double normals*, i.e., pairs $\{x, y\} \subset X$ such that the (open) closed parallel slab with boundary planes orthogonal to the segment $[x, y]$, and containing $x$ or $y$, respectively, contains $X \setminus \{x, y\}$. In the plane their results are similar to those about $a(X)$ and $sa(X)$, but in higher dimensions they are not.

Erdős and Füredi [11] constructed an *acute set* of cardinality $\lfloor \text{const} \cdot (2/\sqrt{3})^d \rfloor$ $= \lfloor \text{const} \cdot 1, 1547 \ldots^d \rfloor$ in $\mathbb{R}^d$ (i.e., all angles determined by the points of such a set are acute, that is stronger than the property of *strict antipodality*, i.e., that each pair of them is strictly antipodal), and they announced that this cardinality can be increased to $\lfloor (2^{1/4} - o(1))^d \rfloor = \lfloor (1.1892 \ldots - o(1))^d \rfloor$. A variant of this was used in [20]. The Erdős-Füredi Theorem was a surprise at that time, since before that only a linear lower bound had been known for the maximal cardinality of a strictly antipodal set in $\mathbb{R}^d$, which was conjectured to be optimal, see [11] or [16, § 7.4]. This was improved by V. Harangi [17], who constructed an acute set of cardinality $\lfloor \text{const} \cdot (144/23)^{d/10} \rfloor = \lfloor \text{const} \cdot 1, 2013 \ldots^d \rfloor$ in $\mathbb{R}^d$. Barvinok-Lee-Novik [1, Theorem 4.1] constructed a still larger strictly antipodal set in $\mathbb{R}^d$, namely one with cardinality $3^{\lfloor d/2-1 \rfloor} - 1 \geq \text{const} \cdot (\sqrt{3})^d$, for $d \geq 4$. D. Zakharov [37] recently constructed an acute set in $\mathbb{R}^d$, of cardinality at least the $d$-th Fibonacci number $\geq [(1 + \sqrt{5})/2]^d$. All these were superceded by an Anonymous from Ukraine for $d = 4, 5$, unpublished, and B. Gerencsér-Harangi [13] for all $d$, who recently constructed an acute – hence strictly antipodal – set in $\mathbb{R}^d$, of cardinality $2^{d-1} + 1$ (cf. also [14]). Thus they disproved the conjecture of P. Erdős ([11]) that this cardinality would be at most $(2 - \varepsilon)^d$, for some $\varepsilon > 0$. It is interesting to observe that for $1 \leq d \leq 3$, where the maximal cardinality of a strictly antipodal set is known, this maximum is just $2^{d-1} + 1$. A variant of this question was investigated by Kupavskii-Zakharov [19]: If the maximal angle between points of $X \subset \mathbb{R}^d$ is at most $\alpha$ $(< \pi/2)$, where $\alpha$ is fixed and $d \to \infty$, then how large can $|X|^{1/d}$ be? They constructed an $X$ with $|X| = \left(\sqrt{2} + o(1)\right)^d$ for some $\alpha$, but it is not clear if this is optimal or not.

For the weaker property of strict antipodality, [6, Lemma 9.11.2] (due to the last mentioned author of this paper) constructed a strictly antipodal set in $\mathbb{R}^d$, of cardinality $3^{\lfloor d/3 \rfloor} \geq \text{const} \cdot (3^{1/3})^d = \text{const} \cdot 1.4422...^d$, and [6, p. 271] announced the result, due to the same author, that here $3^{1/3}$ can be replaced by $5^{1/4} = 1.4953...$. (In our paper we will point out the method for obtaining this stronger estimate, and in fact we will obtain even $4^{1/3} = 1.5874...$, rather than $5^{1/4}$.)

Erdős and Füredi [11] proved also a stronger theorem. What happens if we require that the angles determined by the points should have a smaller upper bound? Clearly, if all angles are at most $\pi/3$, then all angles are exactly $\pi/3$, and we have the vertices of a regular simplex of dimension at most $d$. Hence the number of points is at most $d+1$. However, Erdős-Füredi [11] proved with the probabilistic method that for any $\varepsilon > 0$ there is an $f(\varepsilon) > 1$, such that there exists a set $X_\varepsilon \subset \mathbb{R}^d$, of cardinality at least $f(\varepsilon)^d$, where each angle determined by $X$ is at most $\pi/3 + \varepsilon$.

More exactly, [11] gave an example of $\Omega[(1+\text{const} \cdot \varepsilon^2)^d]$ points in $\mathbb{R}^d$, with quotient of maximal and minimal distances at most $1 + \varepsilon + o(\varepsilon)$. The optimal constant from their method was determined by [20], and was found to be $4/e^2 + o(1) = 0.5413... + o(1)$. Later Frankl and Maehara [12] gave a simple geometric proof for the same statement (just via the greedy algorithm), however, with a slightly worse constant, namely $1/2$.

More general concepts were introduced by Csikós-Kiss-Swanepoel-de Wet in [8] (cf. also Swanepoel-Valtr [35]). They called them antipodality and strict antipodality, but to distinguish them from our concepts, we call them weak antipodality and weak strict antipodality, respectively. Let $\{X_1, \ldots, X_n\}$ be a family of sets in $\mathbb{R}^d$ with $\text{aff}\left(\bigcup_{i=1}^{n} X_i\right) = \mathbb{R}^d$. This family is *weakly antipodal*, or *weakly strictly antipodal*, if for any $i \neq j$, $1 \leq i, j \leq n$, and any $x_i \in X_i, x_j \in X_j$, we have that $x_i, x_j$ are antipodal, or strictly antipodal, with respect to the set $X := \bigcup_{i=1}^{n} X_i$.

[20, Theorem 2] constructed a set $X \in \mathcal{C}_3^s(n)$ with $sa(X) = \lfloor n^2/3 \rfloor$, disproving a conjecture of Grünbaum (see [15] and [16]), and conjectured that the upper bound is $n^2/3 + O(1)$ (p. 461, second part of the conjecture). This was refuted by Csikós-Kiss-Swanepoel-de Wet [8, Theorem 3], which gave an analogous example with $\lfloor (n^2 + n - 2)/3 \rfloor$ points. (Namely, they gave in $\mathbb{R}^3$ a weakly strictly antipodal family of four sets, one being a singleton, and the other three ones being non-trivial $C^1$ arcs. Place on each of the three arcs $\lfloor (n-1)/3 \rfloor$ or $\lceil (n-1)/3 \rceil$ points, altogether $n-1$ points, and the $n$'th point should be the singleton set. Then we obtain a set $X \in \mathcal{C}_3^s(n)$, Then $sa(X) = \lfloor (n^2 + n - 2)/3 \rfloor$.) However, the following conjecture is still feasible:

**Conjecture I.** $\max\{sa(X) \mid X \in \mathcal{A}_3(n), \text{ or } X \in \mathcal{C}_3^s(n)\} = n^2/3 + O(n)$.

For $\mathbb{R}^d$, with $d \geq 4$, [20, Theorem 3] constructed a set $X$ of $n$ points in $\mathbb{R}^d$, with aff $X = \mathbb{R}^d$ and $sa(X) \geq (1 - \text{const}/1.0044\ldots^d)n^2/2 - O(1)$. This was very much improved by Barvinok-Lee-Novik [1, Theorem 4.1], who constructed in $\mathbb{R}^d$, where $d \geq 4$, a set $X$ of $n$ points with aff $X = \mathbb{R}^d$ and $sa(X) \geq \left(1 - 1/(3^{\lfloor d/2-1 \rfloor} - 1)\right) n^2/2 - O(1) \geq \left(1 - \text{const}/(\sqrt{3})^d\right) n^2/2 - O(1)$. For further results we refer to the survey [22].

**3.** Now we turn to citing results that are directly connected to the subject of this paper.

**A.** Brass [7, Theorem 2] proved

$$\max\{a(X) \mid X \in \mathcal{C}_2(n)\} = \lfloor n^2/4 \rfloor + 2 \tag{1}$$

(in the equivalent form that this inequality holds for any norm on the plane for the number of diametral pairs; and thus refuting the first part of the conjecture from [20, p. 461], namely $(n^2/2)(1-1/2^{d-1})$ for $\mathbb{R}^d$). [20, Theorem 2] proved

$$\begin{aligned} \lfloor 3n^2/8 \rfloor + 4 &\leq \max\{a(X) \mid X \in \mathcal{C}_3(n)\} \\ &\leq 7n^2/16 \text{ (the lower bound holding for } n \geq 8\text{)}. \end{aligned} \tag{2}$$

Observe that, for $d + 1 \leq n \leq 2^d$, we have $\max\{a(X) \mid X \subset \mathbb{R}^d, \text{aff } X = \mathbb{R}^d, |X| = n\} = n(n-1)/2$, since we can choose $n$ vertices of a cube. Therefore, we may suppose $n > 2^d$. As a generalization of (1), we pose

**Conjecture II.** Let $X \in \mathcal{A}_d(n)$ or $X \in \mathcal{C}_d(n)$, and let $n > 2^d$. Then $\max a(X)$ is attained in the following case. We take a cube, and on each of $2^{d-1}$ parallel edges we choose $\lfloor n/2^{d-1} \rfloor$ or $\lceil n/2^{d-1} \rceil$ points, altogether $n$ points, which include both endpoints of the respective edge.

Observe that points $x_i, x_j \in X$ on different parallel edges are antipodal, but also the two endpoints of any edge form an antipodal pair. For $2^{d-1} | n$ we have thus $(n^2/2)(1 - 1/2^{d-1}) + 2^{d-1}$ antipodal pairs (and for $2^{d-1} \nmid n$ the first summand has to be replaced by the corresponding Turán number). We remark that an upper bound for $\max a(X)$, close to this conjecture, namely $(n^2/2)(1 - 1/2^d)$, for $2^d | n$, was given in [20, Theorem 3, III] (again, for $2^d \nmid n$, this has to be replaced by the corresponding Turán number). Moreover, as follows from [8, Theorem 1], for $d = 3$ and $k = 1$ this conjecture is true; cf. also our Corollary 4.2.

Nguyễn-V. Soltan [25, Theorem 2] proved in particular

$$n \geq d+1, \ X \in \mathcal{C}_d^s(n) \implies a(X) \geq d(d+1)/2, \text{ with} \\ \text{equality if and only if } n = d+1, \text{ and then } P \text{ is a simplex.} \tag{3}$$

V. Soltan-Nguyễn [33] and Nguyễn-V. Soltan [26, Theorems 1 and 2] proved a generalization of (3), namely

$$\underline{a}_d^s(n) = n + d(d+1)/2 - 1, \\ \text{with equality for } d = 2 \text{ if and only if } P \text{ has no} \\ \text{parallel sides, and for } d \geq 3 \text{ if and only if } P \text{ has } d+1 \tag{4} \\ \text{pairwise antipodal vertices such that any further vertex} \\ \text{is antipodal to a unique one from these } d+1 \text{ vertices.}$$

Thus equality holds, e.g., if $X$ consists of the $d+1$ vertices of a simplex, and still $n - d - 1$ points, each very close to the barycentre of some facet (depending on the point), all points being in strictly convex position.

**B.** The result

$$\underline{sa}_d^s(n) = \lceil n/2 \rceil \text{ for all } n \geq 2d \text{ even, and all } n \geq 4d - 1 \text{ odd} \tag{5}$$

was obtained by Nguyễn [24], deposited but unpublished. Later the case of even $n$ in (5) was obtained in a trivial way in [20, p. 458], and the case of odd $n$ in (5) was obtained in [21, p. 188] in an easy way. Nguyễn-V. Soltan [25, Theorem 3] proved

$$n \geq d+1 \text{ and } X \in \mathcal{C}_d^s(n) \implies sa(X) \geq d, \text{ with equality} \\ \text{if and only if } n = 2d \text{ and } P \text{ is a cross-polytope.} \tag{6}$$

Nguyễn [24] and Nguyễn-V. Soltan [26, Theorem 3] showed that

$$X \in \mathcal{C}_d^s(n) \implies sa(X) \geq \lceil n/2 \rceil. \\ \text{For } n \text{ even, equality holds if and only if} \\ n \geq 2d, \text{ and vert } P \text{ is the union of } n/2 \text{ disjoint pairs,} \tag{7} \\ \text{such that for any pair the respective support cones are} \\ \text{centrally symmetric images of each other with respect to } 0.$$

Moreover, they clarified the cases of equality also for odd $n$, which holds if and only if $n \geq 4d - 1$, and a slightly more difficult geometric condition holds. As Nguyễn and V. Soltan (see [33] and [26]) observed, this solves the question of $\underline{sa}_2^s(n)$, namely

$$\underline{sa}_2^s(3) = 3 \text{ and, for } n \geq 4, \ \underline{sa}_2^s(n) = \lceil n/2 \rceil \tag{8}$$

(cf. also the explanation in [20, § 1, p. 459]).

The question of $\underline{sa}_3^s(n)$ was solved by Nguyễn [24], deposited but unpublished, and later a proof was published by [21, § 3]. The result is:

$$\begin{aligned} &\underline{sa}_3^s(4) = \underline{sa}_3^s(5) = \underline{sa}_3^s(7) = \underline{sa}_3^s(9) = 6, \\ &\underline{sa}_3^s(n) = \lceil n/2 \rceil \text{ for all other } n\text{'s} . \end{aligned} \tag{9}$$

Here $\underline{sa}_3^s(4) = 6$ is trivial, the proof of $\underline{sa}_3^s(5) = 6$ was easy, but $\underline{sa}_3^s(7) \geq 6$ and $\underline{sa}_3^s(9) \geq 6$ were proved in a longer way. However, [21], §1, p. 186 observed that these inequalities also follow from a result of Yost [36], that just appeared that time. Then there remain to show $\underline{sa}_3^s(5), \underline{sa}_3^s(7), \underline{sa}_3^s(9) \leq 6$, that is easily done by some simple examples, cf. [21, p. 192].

The result of Yost [36] was the following. We call a centrally symmetric convex polytope $Q \subset \mathbb{R}^d$ *irreducible*, if its 0-symmetric translate is not of the form $(P - P)/2$, where $P \subset \mathbb{R}^d$ is a convex polytope that is no translate of $Q$. *Reducible* means: not irreducible. Yost [36, Theorem 11] states that

$$\begin{aligned} &\text{for } d \geq 3, \text{ any centrally symmetric polytope} \\ &\quad \text{with less than } 4d \text{ vertices is irreducible.} \end{aligned} \tag{10}$$

For further discussion of this property cf. § 3 of this paper, where we use it for the determination of $\underline{sa}_d^s(n)$, for $2d + 1 \leq n \leq 4d - 1$ odd.

Therefore, the cases not settled by (5) are: $d + 2 \leq n \leq 2d - 1$, and $n$ odd with $2d + 1 \leq n \leq 4d - 3$. In our paper we solve the second case, and for the first case we determine $\underline{sa}_d^s(2d - 1)$, and we give lower and upper bounds to $\underline{sa}_d^s(n)$, which will be of the same order of magnitude, for $d + 2 \leq n \leq (2 - \varepsilon)d$.

Our results below about $\underline{sa}_d^s(n)$ will fill in certain gaps in [20] and [21]. Namely in these papers, in pages 459, and 186, respectively, there were announced several statements without proofs, namely $(2_2)$ to $(2_4)$, about $\underline{sa}_d^s(n)$ for the cases not settled by (5). In this paper all these announced results will be proved. So there remains the only

**Problem.** Determine $\underline{sa}_d^s(n)$ for $d+2 \leq n \leq 2d-2$ exactly.

**C.** [21] considered a generalization of the concept of antipodal pairs of points. Let $S^k = \{s_1^k, \ldots, s_n^k\}$ be a finite set of $k$-simplices in $\mathbb{R}^d$, with $\mathrm{aff}\left(\bigcup_{i=1}^n s_i^k\right) = \mathbb{R}^d$. We call $S^k$ *$k$-antipodal, or strictly $k$-antipodal*, if for any $i \neq j$, where $1 \leq i, j \leq n$, there are different parallel supporting hyperplanes $H', H''$ of $P^k := \mathrm{conv}\left(\bigcup_{i=1}^n s_i^k\right)$ such that $s_i^k \subset H'$ and $s_j^k \subset H''$, or $s_i^k = P^k \cap H'$ and $s_j^k = P^k \cap H''$, respectively. [21, p. 187] formulated the conjecture that a $k$-antipodal set of $k$-simplices has at most $2^{d-k}$ elements. This number is attained, e.g., for $2^{d-k}$ $k$-simplices on $2^{d-k}$ parallel $k$-faces of a cube. Observe that the case $k = 0$ of this conjecture is the Danzer-Grünbaum Theorem [9]. Moreover, for the case $k = d-1$ we have an evident positive answer, so we may suppose $1 \leq k \leq d-2$. Formerly, also I. Bárány and V. Soltan had formulated this conjecture.

Observe that this number $2^{d-k}$ is attained not only in the above written case. Namely, if in $\mathbb{R}^d$ we have some two sets of $2^{d-k}$ $k$-simplices that are antipodal, using the cube $[-1, 1]^d$, then in $\mathbb{R}^{d+1}$ we can embed these examples in $[-1, 1]^d \times \{-1\}$, and $[-1, 1]^d \times \{1\}$, even in any rotated position. Thus, in $\mathbb{R}^3$ we have two combinatorially different examples, and the number of examples becomes at least squared at each step. Thus the number of combinatorially different examples is at least $2^{2^{d-3}}$. The symmetry group of the cube has $2^d \cdot d!$ elements. So the number of examples being not combinatorially isomorphic is still at least $2^{2^{d-2}}/[2^d \cdot d!]$. This makes this question complicated. Therefore, we pose the simpler

**Conjecture III.** For $1 \leq k \leq d-2$, the number of pairwise $k$-antipodal affine $k$-subspaces in $\mathbb{R}^d$ is at most $2^{d-k}$, with equality if and only if they are the affine hulls of some $2^{d-k}$ parallel $k$-faces of a parallelepiped.

Here $k$-antipodality of affine $k$-subspaces is defined as for $k$-simplices.

By B. Gerencsér-Harangi [13], there exists an acute, hence strictly antipodal set $S \subset \mathbb{R}^{d-k}$ with $|S| = 2^{d-k-1} + 1$. Then the translates of the orthocomplement of $\mathbb{R}^{d-k}$ in $\mathbb{R}^d$, containing the points of $S$, form a strictly antipodal set of affine $k$-subspaces, of cardinality $2^{d-k-1} + 1$, which is quite close to the value in Conjecture III. Of course, this example yields strictly antipodal sets of $k$-simplices in $\mathbb{R}^d$, of the same cardinality, as well.

For the case $d = 3$ and $k = 1$ a partial result was proved in [21, Proposition 2]:

> If $S^1 = \{s_1^1, \ldots, s_n^1\}$ is a strictly antipodal set of segments in $\mathbb{R}^3$ which are pairwise skew, then $n \leq 3$, and this bound is sharp. \hfill (11)

In § 3 we will show that this implies that an antipodal, or strictly antipodal, set of segments in $\mathbb{R}^3$ consists of at most four, or three segments, respectively, and we describe also the cases of equality. This inequality, for the strictly antipodal case, is given, in a more general form, in Csikós-Kiss-Swanepoel-de Wet [8, Theorem 2], as follows.

> In $\mathbb{R}^3$ there is no weakly strictly antipodal set consisting of four $C^1$ arcs. \hfill (12)

However, the cases of equality (for three arcs) are not discussed there.

Observe that four segments on four parallel edges of a parallelepiped in $\mathbb{R}^3$ are an antipodal set of segments in $\mathbb{R}^3$. The same statement holds for any subsegments of two parallel edges of a face of the parallelepiped in $\mathbb{R}^3$, and of two parallel edges of the opposite face of our parallelepiped, which are not parallel to the first two edges. Csikós-Kiss-Swanepoel-de Wet [8, Theorem 1] assert that

> any weakly antipodal family of four segments in $\mathbb{R}^3$ is of one of the above described two types. \hfill (13)

Theorem 4 of [8] asserts that, for $\mathbb{R}^3$,

> for $n = 6$, and $\min_{1 \leq i \leq 6} |X_i|$ sufficiently large, $\{X_1, \ldots, X_6\}$ cannot be weakly antipodal. \hfill (14)

Thus a weakly antipodal set of segments consists of at most five segments. However, (13) easily implies that

> a weakly antipodal set of segments in $\mathbb{R}^3$ consists of at most four segments. \hfill (15)

In fact, suppose that we have five such segments. It will suffice to consider only the directions of the segments. If all are parallel, we project them along their direction to a plane, obtaining an antipodal set of five points, contradicting the Danzer-Grünbaum theorem. If not all segments are parallel,

$$\text{any four of them should be parallel, or should form two parallel pairs not parallel to each other.} \quad (16)$$

Therefore, there are no three pairwise not parallel segments. Thus the segments have two directions, and the numbers of segments of the same directions should be 4 and 1, or 3 and 2, both contradicting (16).

Thus (13) is stronger than our Theorem 4, antipodal case, but it seems that the deduction of our Theorem 4, antipodal case, is simpler than that of Theorem 1 of [8].

For some other aspects of antipodality, we refer to [30], [3] and [4].

**4.** Now we list all the possible questions that are of the type investigated in this paper, and make some remarks on them. Observe that we considered $a(X), sa(X)$, their minima and maxima, and $X$ was in $\mathcal{A}_d(n)$, $\mathcal{C}_d(n)$, or $\mathcal{C}_d^s(n)$. Therefore, we have twelve questions, which we list below. We write $P := \operatorname{conv} X$.

**A)** Question of $\min a(X)$, for $X \in \mathcal{A}_d(n)$, $X \in \mathcal{C}_d(n)$, or $X \in \mathcal{C}_d^s(n)$. For $X \in \mathcal{A}_d(n)$, we can have vertices of a simplex, and $n - d - 1$ points in its interior, giving $a(X) = (d+1)d/2$. On the other hand, $a(X) = a(X \cap \operatorname{bd} P) \geq a(\operatorname{vert} P) \geq (d+1)d/2$ by [33] (cf. also [26]). For $X \in \mathcal{C}_d(n)$ cf. our Theorem 1.2. For $X \in \mathcal{C}_d^s(n)$, the answer is $n + d(d-1)/2 - 1$, cf. [26].

**B)** Question of $\max a(X)$, for $X \in \mathcal{A}_d(n)$, $X \in \mathcal{C}_d(n)$, or $X \in \mathcal{C}_d^s(n)$. Since $a(X) = a(X \cap \operatorname{bd} P)$, the answer to the first of these questions is $\max\{\overline{a}_d(m) \mid d+1 \leq m \leq n\}$. The second of these questions see in [20] and in Conjecture II in this paper. (If Conjecture II is true, then the answers to the first and second questions coincide.) The third of these questions was treated in [20].

**C)** Question of $\min sa(X)$, for $X \in \mathcal{A}_d(n)$, $X \in \mathcal{C}_d(n)$, or $X \in \mathcal{C}_d^s(n)$. Since $sa(X) = sa(\operatorname{vert} P)$, the first two of these quantities equal $\min\{\underline{sa}_d^s(m) \mid m \leq n\}$, and the third one is $\underline{sa}_d^s(n)$.

**D)** Question of $\max sa(X)$, for $X \in \mathcal{A}_d(n)$, $X \in \mathcal{C}_d(n)$, or $X \in \mathcal{C}_d^s(n)$. Here we have $sa(X) = sa(\operatorname{vert} P)$, hence the first and second of these quantities equal $\max \overline{sa}_d^s(m) \mid d+1 \leq m \leq n\}$. However, this equals $\overline{sa}_d^s(n)$, since this quantity is (strictly) increasing with $n$. (In fact, to any $X \in \mathcal{C}_d^s(n)$, we can add a further point $y$ outside $P$, but very close to the barycentre of a facet. Then strictly antipodal pairs of $X$ will remain strictly antipodal in $X \cup \{y\}$, and there is an additional strictly antipodal pair in $X \cup \{y\}$, containing $y$.)

## 2 Theorems

We improve the upper estimate in (2), i.e., in [20], Theorem 2.

**Theorem 1.1.** *We have*

$$\max\{a(X) \mid X \subset \mathbb{R}^3,\ \mathrm{aff}\, X = \mathbb{R}^3,\ X \text{ is in convex position},\ |X| = n\} \leq 2n^2/5 + O(n^c),\ \text{for some } c < 2.$$

**Theorem 1.2.** *Let $d \geq 2$, and $n \geq d+1$ be integers. Then*

$$\min\{a(X) \mid X \subset \mathbb{R}^d,\ \mathrm{aff}\, X = \mathbb{R}^d,\ X \text{ is in convex position},\ |X| = n\} = n + d(d-1)/2 - 1\,.$$

The following observation will be used several times.

**Lemma 2.1.** *Let $Y \subset \mathbb{R}^d$ be a finite set, in convex position, or strictly convex position, respectively, whose affine hull is $(d-1)$-dimensional, and let $x \in \mathbb{R}^d \setminus (\mathrm{aff}\, Y)$, and $X = Y \cup \{x\}$. Then, considering $a(Y)$ and $sa(Y)$ in $\mathrm{aff}\, Y$, we have*

$$a(X) = a(Y) + |Y| \ \text{ or } \ sa(X) = sa(Y) + |Y|,\ \text{respectively.}$$

**Definition.** Let $d \geq 3$, and $n \geq d+1$ be integers. We let

$$(sa_d^s)'(n) := sa_d^s(n) \text{ for } n \neq 2d,\ \text{and } (sa_d^s)'(2d) := \min\{sa(X) \mid X \in \mathcal{C}_d^s(2d),\ \text{and conv } X \text{ is not a cross-polytope}\}.$$

We give the simple

**Lemma 2.2.** *Let $X \subset \mathbb{R}^d$ be a finite set with $\mathrm{aff}\, X = \mathbb{R}^d$, and let $\pi : \mathbb{R}^d \to \mathbb{R}^{d-1}$ be the projection $(x_1, \ldots, x_{d-1}, x_d) \mapsto (x_1, \ldots, x_{d-1})$. Let $Y = \pi(X)$. If $y_1, y_2 \in Y$ are strictly antipodal (in $\mathrm{aff}\, Y = \mathbb{R}^{d-1}$), then the highest point* in $X \cap \pi^{-1}(y_1)$ and the lowest point of $X \cap \pi^{-1}(y_2)$ form together a strictly antipodal pair in $X$. However, if $|X \cap \pi^{-1}(y_1)|, |X \cap \pi^{-1}(y_2)| \geq 2$, then the highest points (lowest points) of $X \cap \pi^{-1}(y_1)$ and $X \cap \pi^{-1}(y_2)$ do not form a strictly antipodal pair in $X$.

We remark that [21] already tacitly used the first statement of this lemma, without explicitly mentioning it.

**Theorem 2.3.** *For $\underline{sa}_d^s(d+k)$, where $2 \leq k \leq d-1$, we have the following upper bounds.*

(1) *For $2 \leq k \leq \lfloor (2d-2)/3 \rfloor$, we have $\underline{sa}_d^s(d+k) \leq (d(d+1)/2)$ $+k-1 = \Theta(d^2)$,*

(2) *for $\lceil (2d-2)/3 \rceil \leq k \leq d-1$, we have $\underline{sa}_d^s(d+k) \leq k + (d+k) \times (d+k-1)/2 - 2k(2k-1)/2 = k + (d-k)(d+3k-1)/2 = \Theta((d-k)d)$.*

**Remark.** Both examples (namely in (1) and (2)) work for all $k \in \{2, \ldots, d-1\}$. For $1 < k < (2d-2)/3$ the formula in the first example is better, for $k = (2d-2)/3$ (whether this is an integer, or not) the formulas in the two examples give the same value, and for $(2d-2)/3 < k \leq d-1$ the formula in the second example is better, as readily checked. In case (1), the upper bound is $d^2/2 + O(d)$. In case (2), with $k = \lfloor xd \rfloor$ for $x \in [2/3, 1)$ fixed, the upper bound is $(d^2/2)f(x) + O(d)$, where $f(x) := 3(1-x)(x+1/3)$. Here $f(0) = f(2/3) = 1$ and $f(1) = 0$, and $f(x)$ is strictly decreasing for $x \in [2/3, 1]$.

The second statement of the following theorem will be sharpened in Theorem 2.15. However, the proof of this sharpening will use itself the second statement of Theorem 2.4.

**Theorem 2.4.** *Let $d \geq 3$, and let $n \in [2d+1, 4d-1]$ be an odd integer. Then*

$$\underline{sa}_d^s(n) = 2d.$$

*Let $d \geq 3$, and let $d+2 \leq n \leq 2d-1$. Then*

$$\underline{sa}_d^s(n) \geq 2d.$$

Let $d \geq 3$, and let $X \in \mathcal{C}_d^s(2d)$. Then

$$\text{either } P = \text{conv} X \text{ is a cross-polytope,}$$
$$\text{and then } sa(X) = d, \text{ or else } sa(X) \geq 2d.$$

In the following proposition we show that the minimum $\underline{sa}_3^s(5) = 6$ is attained only for a pyramid over a paralellogram. This sharpens the statement $\underline{sa}_3^s(5) = 6$ in [21, p. 188, 2A], cf. also (9).

**Proposition 2.5.** *We have $\underline{sa}_3^s(5) = 6$. For $X \in \mathcal{C}_3^s(5)$ we have $sa(X) = 6$ if and only if*
*1) $P = \text{conv} X$ is a pyramid over a parallelogram.*
*For all other cases we have:*
*2) if $P$ is a pyramid over a trapezoid that is no parallelogram, then $sa(X) = 7$,*
*3) if $P$ is a pyramid over a quadrangle that is no trapezoid, then $sa(X) = 8$,*
*4) if $P$ is a bipyramid over a triangle, then $7 \leq sa(X) \leq 10$, with both bounds attained.*

We know $\underline{sa}_3^s(6) = 3$ (see [20, p. 458, (1)], [21, p. 186, (1)] and (9)). The following proposition, which is the special case $d = 3$ of Theorem 2.4, will sharpen this statement.

**Proposition 2.6.** *We have $\underline{sa}_3^s(6) = 3$. For $X \in \mathcal{C}_3^s(6)$ we have $sa(X) = 3$ if and only if $P = \text{conv} X$ is a cross-polytope. Else we have $sa(X) \geq 6$.*

We have the following lower estimates.

**Theorem 2.7.** *Let $d \geq 3$ and $k \in \{2, \ldots, d-1\}$ be integers. Let $X \in \mathcal{C}_d^s(d+k)$. Let*

$$sa(X) \geq c \geq \min\{(\underline{sa}_{d-1}^s)'(d-1+k') \mid 1 \leq k' \leq k\}.$$

*Then we have the sharper estimate*

$$sa(X) \geq \min\{(\underline{sa}_{d-1}^s)'(d-1+k') \mid 1 \leq k' \leq k\} + \lceil 2c/(d+k) \rceil,$$

*unless we have*

$$k = d - 1, \text{ and } X \text{ is a pyramid over a } (d-1)\text{-dimensional} \quad (17)$$
$$\text{cross-polytope, and then } sa(X) = 3(d-1).$$

*In particular, we have*

$$sa(X) \geq \lceil \min\{(\underline{sa}^s_{d-1})'(d-1+k') \mid 1 \leq k' \leq k\} \times (d+k)/(d+k-2) \rceil,$$

*unless we have (17).*

With the help of Theorem 2.7, we strengthen the second statement of Theorem 2.4, for $d \geq 4$.

**Proposition 2.8.** *Let $d \geq 4$, and $k \in \{1, \ldots, d-1\}$ be integers. Then*

$$\underline{sa}^s_d(d+k) \geq 2d+1.$$

The following theorem will finish the determination of $\underline{sa}^s_d(n)$ for $d = 4$ and all $n \geq d+1$, except for $\underline{sa}^s_4(7) = 9$, which will follow from Theorem 2.15. Observe that (5) settles the cases $n \geq 2d$ even and $n \geq 4d+1$ odd, while Theorem 2.4 settles the case $n \in [2d+1, 4d-1]$ odd. So there remain the cases $d+2 \leq n \leq 2d-1$ (the case $n = d+1$ being trivial). For $d = 4$ this means that there remain the cases $6 \leq n \leq 7$.

**Theorem 2.9.** *We have*

$$\underline{sa}^s_4(6) = 11.$$

*The minimum $\underline{sa}^s_4(6) = 11$ is attained, e.g., for a two-fold (i.e., twice iterated) pyramid over a parallelogram, and also for the convex hull of a simplex and an additional vertex very close to the barycentre of one of its facets.*

**Proposition 2.10.** *Let $d \geq 2$ and $k \geq 1$ be integers, and let $X \in \mathcal{C}^s_d(d+k)$. Suppose that in the graph $G$ with vertex set $X$ and edge set the strictly antipodal pairs of $X$, the maximal degree is $\ell$. Then there exists an orthogonal projection $\pi$ from $\mathbb{R}^d$ to a linear $(d-1)$-subspace of itself, such that for $P := \operatorname{conv} X$, $Q := \pi P$ and $Y := \operatorname{vert} Q \in \mathcal{C}^s_{d-1}(|Y|)$ we have*
*(1) $|Y| \leq |X| - 1$, and*
*(2) $sa(X) \geq sa(Y) + \ell \geq \min\{\underline{sa}^s_{d-1}(d-1+k') \mid 1 \leq k' \leq k\} + \ell.$*

**Theorem 2.11.** *Suppose that the hypotheses of Proposition 2.10 hold. Then*

$$sa(X) \geq \lceil \min\{\underline{sa}^s_{d-1}(d-1+k') \mid 1 \leq k' \leq k\} \cdot (d+k)/(d+k-2) \rceil.$$

*Hence*

$$\underline{sa}_d^s(d+k) \geq \lceil \min\{\underline{sa}_{d-1}^s(d-1+k') \mid 1 \leq k' \leq k\} \times \\ (d+k)/(d+k-2)\rceil. \tag{18}$$

*Equivalently,*

$$\min\{\underline{sa}_d^s(d+k') \mid 1 \leq k' \leq k\} \geq \\ \lceil \min\{\underline{sa}_{d-1}^s(d-1+k') \mid 1 \leq k' \leq k\} \cdot (d+k)/(d+k-2)\rceil.$$

**Lemma 2.12.** *If a convex polytope $P \subset \mathbb{R}^d$ contains a unique segment $[a,b]$ of maximal length of some direction, where both $a,b$ are vertices of $P$, then $a$ and $b$ are strictly antipodal in $P$.*

This is true only for both $a, b$ vertices: e.g., for $d = 2$ and $P = \operatorname{conv}\{a = (0,0), (0,2), (2,0)\}$ and $b = (1,1)$ this is false.

**Proposition 2.13.** *Let $d \geq 2$ and $k \geq 1$ be integers, and let $X \in \mathcal{C}_d^s(d+k)$. Suppose that in the graph $G$ with vertex set $X$ and edge set the strictly antipodal pairs of $X$, the maximal degree is $\ell \leq d+k-2$. Then there exists an orthogonal projection $\pi$ from $\mathbb{R}^d$ to a linear $(d-1)$-subspace of itself, such that for $P := \operatorname{conv} X$, $Q := \pi P$ and $Y := \operatorname{vert} Q$ we have*
*(1) $|Y| \leq |X| - 2$, and*
*(2) $sa(X) \geq sa(Y) + \ell \geq \min\{\underline{sa}_{d-1}^s(d-1+k') \mid 1 \leq k' \leq k-1\} + \ell$.*

**Proposition 2.14.** *Suppose that the hypotheses of Proposition 2.13 hold. Then*

$$sa(X) \geq \lceil \min\{\underline{sa}_{d-1}^s(d-1+k') \mid 1 \leq k' \leq k-1\} \times \\ (d+k)/(d+k-2)\rceil.$$

**Theorem 2.15.** *Let $d \geq 2$ and $1 \leq k \leq d-1$ be integers. Then*

$$\underline{sa}_d^s(d+k) \geq 3(d-1).$$

*In particular,*

$$\underline{sa}_d^s(2d-1) = 3(d-1).$$

**Remark.** It would be tempting to use Theorem 2.15 as base of induction, and the last formula in Theorem 2.11 as induction step, to get a lower estimate of $\underline{sa}_d^s(d+k)$, inductively for all $d \geq k+1$. However, the formula thus obtained, namely $\lceil 3(d+k-1)(d+k)/(2(2k+1)) \rceil = \Theta(d^2/k)$, is never better than the lower estimate in the later Theorem 2.18, as can be checked by a simple calculation. Therefore we will not get into it.

Since $\underline{sa}_d^s(d+1) = d(d+1)/2$, and $\underline{sa}_d^s(2d-1) = 3(d-1)$ by Theorem 2.15, it suffices to investigate $\underline{sa}_d^s(d+k-1)$ only for $2 \leq k \leq d-2$.

**Theorem 2.16.** *Let $d \geq 4$ and $2 \leq k \leq d-2$ be integers. Then*

$$\underline{sa}_d^s(d+k) \geq \min \begin{cases} \min\{\underline{sa}_{d-1}^s(d-1+k') \mid 1 \leq k' \leq k\} + (d+k-1), \\ \min\{\lceil \underline{sa}_{d-1}^s(d-1+k') \rceil \mid 1 \leq k' \leq k-1\} \cdot (d+k)/(d+k-2) \end{cases}. \tag{19}$$

**Lemma 2.17.** *Let $d \geq 4$ and $2 \leq k \leq d-2$ be integers. Let*

$$f(d,k) := (d^2 - k^2 + d + k)/2 = (d+k)(d-k+1)/2$$

*(which is integer valued, and is a decreasing function of $k$). Then*

$$f(d-1,k) + (d+k-1) \geq f(d-1,k-1) \cdot (d+k)/(d+k-2).$$

**Theorem 2.18.** *Let $d \geq 2$ and $k \in \{1, \ldots, d-1\}$ be integers. Then*

$$\underline{sa}_d^s(d+k) \geq (d^2 - k^2 + d + k)/2 = \Theta((d-k)d).$$

**Remark.** Letting $k = \lfloor d \cdot x \rfloor$, with $x \in (0,1)$, the lower bound in Theorem 2.16 is $(d^2/2)(1-x^2) + O(d)$. It is interesting to observe that the proof of Theorem 2.18 only will use (19) from Theorem 2.16, but will not use (18) from Theorem 2.11. However, even using also (18) would not give a better lower estimate. Also using the decreasing property of the function $k \mapsto f(d,k)$ from Lemma 2.17, for this is suffices to show the inequality $f(d,k) \geq f(d-1,k) \cdot (d+k)/(d+k-2)$. Rearranging this, this turns out to be equivalent to $k \geq 1$, hence it holds.

**Corollary 2.19.** *We have*

$$\underline{sa}_5^s(6) = 15, \ \underline{sa}_5^s(7) = 16, \ 15 \leq \underline{sa}_5^s(8) \leq 17, \ \text{and} \ \underline{sa}_5^s(9) = 12.$$

**Remark.** Observe that in the proof of Corollary 2.19, for $k = 2$ (18) will be better, while for $k = 3$ (19) will be better. So neither of (18) and (19) is stronger than the other one.

We can sum up Theorems 2.3 and 2.18 as follows.

**Proposition 2.20.** *We have*

$$\underline{sa}_d^s(d+k) = \Theta\left((d-k)d\right).$$

*The quotient of the upper estimate of $\underline{sa}_d(d+k)$ in Theorem 2.3 and of its lower estimate in Theorem 2.18 is at most* 2, *which value is sharp.*

**Theorem 2.21.** *Let $D \geq 2$. Suppose that we know some positive integer lower estimates $\underline{sa}_d^s(d+k)^*$ for $\underline{sa}_d^s(d+k)$ for all $d < D$ and all $1 \leq k \leq d-1$. (We may still suppose that $\underline{sa}_d^s(d+k)^* = \underline{sa}_d^s(d+k)$ for those $d, k$ for which we know $\underline{sa}_d^s(d+k)$.) Then a valid lower estimate of $\underline{sa}_D^s(D+k)$, for $1 \leq k \leq D-2$ is*

$$\max\{\min \\ \{\min\{\underline{sa}_{D-1}^s(d-1+k')^* \mid 1 \leq k' \leq k\} + (d+k-1), \\ \min\{\lceil \underline{sa}_{D-1}^s(d-1+k')^* \mid 1 \leq k' \leq k-1\} \cdot (d+k)/(d+k-2)\rceil\}, \\ \min\{\lceil \underline{sa}_{D-1}^s(d-1+k')^* \mid 1 \leq k' \leq k\} \cdot (d+k)/(d+k-2)\rceil\}, \quad (20)$$

*while for $k = D-1$ it is $3(D-1)$. Of course, by iterating this, we obtain valid lower bounds $\underline{sa}_d^s(d+k)^*$ for $\underline{sa}_d^s(d+k)$, for all $d \geq D$ and all $1 \leq k \leq d-1$.*

Of course, if $\underline{sa}_d^s(d+k)^*$ is a decreasing function of $k$, then the formula simplifies: in the innermost minima we always can consider only $k' = k - 1$ or $k' = k$. Theorem 2.21 gives only an algorithmic lower bound for $\underline{sa}_d^s(d+k)$, but we do not have a formula for this algorithmic lower bound.

Although Barvinok-Lee-Novik [1], Zakharov [37], and B. Gerencsér-Harangi [13] already constructed better examples, cf. **2** of the introduction of our paper, it may have some interest to give a smaller strictly antipodal set in $\mathbb{R}^d$, using the idea presented in [6, proof of Lemma 9.11.2], due to the last mentioned author of this paper. Actually we will give an improvement of that construction, but after the proof of Theorem 3 we will point out the original construction.

**Theorem 3.** *In $\mathbb{R}^d$ there exists a strictly antipodal set, affinely spanning $\mathbb{R}^d$, of cardinality*

$$4^{\lfloor d/3 \rfloor} \geq \mathrm{const} \cdot (4^{1/3})^d = \mathrm{const} \cdot 1.5874...^d.$$

We improve (11), i.e., [21, Proposition 2]. In the following theorem, we will drop the upper indices 1. The proof of this theorem will essentially rely on (11).

**Theorem 4.1** *Let $S = \{s_1, \ldots, s_n\}$ be a set of segments in $\mathbb{R}^3$, affinely spanning $\mathbb{R}^3$. If $S$ is antipodal, or strictly antipodal, then $n \leq 4$, or $n \leq 3$, respectively. For the antipodal case we have equality if and only if the segments are subsegments of some four parallel edges of a parallelepiped, or they are subsegments of two parallel edges of a face of a parallelepiped, and of two parallel edges of the opposite face of the parallelepiped, which are not parallel to the first two edges. For the strictly antipodal case we have equality if and only if the segments are subsegments of the lateral edges of a triangular prism, or they are subsegments of the relative interiors of three mutually skew edges of a parallelepiped.*

As mentioned in the introduction, the antipodal case is a simple consequence of [8, Theorem 1], but it seems that we have a simpler deduction. The inequality for the strictly antipodal case is implied by [8, Theorem 2], but the case of equality is not discussed there.

**Corollary 4.2** *An antipodal, or a strictly antipodal set of lines in $\mathbb{R}^3$ consists of at most four, or three lines, respectively. The only cases of equality are four lines spanned by four parallel edges of a parallelepiped, or by the three lateral edges of a triangular prism, respectively.*

## 3 Maximal and minimal numbers of antipodal pairs

We denote by $K_{t_1,\ldots,t_l}$ the complete $l$-partite graph, i.e., a graph whose vertex set is partitioned into classes with $t_1, \ldots, t_l$ vertices, and two vertices are connected by an edge if and only if they belong to different partition classes.

**Proof of Theorem 1.1:** We define a graph $G$ whose vertex set is $X$, and two vertices are connected by an edge if and only if they form an antipodal pair. By [8, Theorem 4], $G$ does not contain the complete 6-partite graph $K_{m,m,m,m,m,m}$ for some positive integer $m$. Therefore, by [5, Corollary 4.7] the number of edges of $G$ is at most

$$(4/5) \cdot n^2/2 + \text{const} \cdot n^{2-1/m},$$

that proves the theorem. □

**Proof of Theorem 1.2:** We write $P := \text{conv}\, X$, and $Q := (P - P)/2$. The antipodal pairs of $X$ are the pairs $x_i, x_j \in X$ whose distance, in the norm $\|\cdot\|_Q$, with unit ball $Q$, is maximal, i.e., is 2.

Now observe that by sufficiently small perturbations of the points in $X$ the number $a(X)$ cannot increase. In fact. move each point of $X$ at most by some small distance $\varepsilon > 0$, obtaining $X'$, and let us write $P' := \text{conv}\, X'$, $Q' := (P' - P')/2$. Then the Hausdorff distance of $P$ and $P'$ is $\max\{|h_P(u) - h_{P'}(u)| \mid u \in S^{d-1}\}$, where $h_P(u) = \max\{\langle x, u\rangle \mid x \in X\}$ and $h'_P(u) = \max\{\langle x', u\rangle \mid x' \in X'\}$, $u \in S^{d-1}$, cf. [29]. Hence, for each $u \in S^{d-1}$, we have that $|h_P(u) - h_{P'}(u)| \le \varepsilon$, that implies $|h_Q(u) - h_{Q'}(u)| \le \varepsilon$. That is, the Hausdorff distance of $Q$ and $Q'$ is at most $\varepsilon$.

Now let $x_i, x_j \in X$, with $\|x_i - x_j\|_Q < 2$. Then, for sufficiently small $\varepsilon$, we have $\|x'_i - x'_j\|_{Q'} < 2$, that is, we have

$$a(X') \le a(X).$$

Now let $K \supset P$ be a strictly convex body whose Hausdorff distance to $P$ is very small. From some point of $\text{int}\, P$ we project radially each $x \in X \subset \text{bd}\, P$ to $\text{bd}\, K$, obtaining $X'$, which is in strictly convex position. Then, for each $x \in X$, the corresponding point $x'$ is very close to $x$. Hence

$$a(X) \ge a(X') \ge n + d(d-1)/2 - 1,$$

the second inequality following from (4). □

# 4 Minimal number of strictly antipodal pairs

**Proof of Lemma 2.1:** We have to observe that the antipodal, or strictly antipodal pairs in $X$ are those of $Y$ (considered in $\text{aff}\, Y$), and any pair containing $x$ (by convex position, or strictly convex position of $Y$, respectively). □

**Proof of Lemma 2.2:** Let $H', H''$ be two different parallel supporting hyperplanes of $Y$ (taken in $\mathbb{R}^{d-1}$), such that $Y \cap H' = \{y_1\}$ and $Y \cap H'' = \{y_2\}$. Then we have two different parallel supporting hyperplanes of $X$, namely $\pi^{-1}(H')$ and $\pi^{-1}(H'')$, such that $X \cap \pi^{-1}(H') = X \cap \pi^{-1}(y_1)$ and $X \cap \pi^{-1}(H'') = X \cap \pi^{-1}(y_2)$. Let, e.g., $y_1 = (1, 0, \ldots, 0), y_2 = (-1, 0, \ldots, 0)$. Let the highest point in $X \cap \pi^{-1}(y_1)$ be $x_1$, and the lowest point of $X \cap \pi^{-1}(y_2)$ be $x_2$. These lie in the $\xi_1 \xi_d$-coordinate plane ($\xi_i$ being the coordinates, the $\xi_1$-axis suitably chosen). Now we perturb $\pi^{-1}(H')$ and $\pi^{-1}(H'')$, such that the relations $x_1 \in \pi^{-1}(H')$ and $x_2 \in \pi^{-1}(H'')$ remain valid, and both the lowest point of $X \cap \pi^{-1}(y_1)$, if different from $x_1$, and the highest point of $X \cap \pi^{-1}(y_2)$, if different from $x_2$, should be strictly between the perturbed hyperplanes. Then all points of $X \setminus (\pi^{-1}(H') \cup \pi^{-1}(H''))$ will also lie strictly between the perturbed hyperplanes.

The second statement is sufficient to be proved for the highest points. Let, for $i = 1, 2$, the highest, and lowest point of $X \cap \pi^{-1}(y_i)$ be $y_i + \lambda_i e_d$, and $y_i + (\lambda_i - c_i)e_d$, where $e_d$ is the standard $d$'th unit vector of $\mathbb{R}^d$, and $c_i > 0$. Let $\varphi$ be a linear functional on $\mathbb{R}^d$, with its unique minimum attained at $y_1 + \lambda_1 e_d$ and unique maximum attained at $y_2 + \lambda_2 e_d$. Then $\varphi(y_1 + \lambda_1 e_d) < \varphi(y_1 + (\lambda_1 - c_1)e_d)$, and $\varphi(y_2 + \lambda_2 e_d) > \varphi(y_2 + (\lambda_2 - c_2)e_d)$. These imply $\varphi e_d < 0$ and $\varphi e_d > 0$, a contradiction. $\square$

**Proof of Theorem 2.3:** The $\Theta(\cdot)$ relation in (2) of the theorem follows from $O(d) + (d-k)\Theta(d) = \Theta((d-k)d)$.

We have to give examples of sets $X$ for which, for $2 \le k \le \lfloor (2d-2)/3 \rfloor$, we have $sa(X) = d(d-1)/2 + d + k - 1$, and for $\lceil (2d-2)/3 \rceil \le k \le d-1$ we have $sa(X) = k + (d+k)(d+k-1)/2 - 2k(2k-1)/2$.

In the first case we recall (4). Here equality holds, e.g., if $X$ consists of the $d+1$ vertices of a simplex, and still $k-1$ points, each very close to the barycentre of some facet, all points being in strictly convex position. Clearly, for this example we have also

$$sa(X) = d(d+1)/2 + k - 1\,.$$

Observe that this example works for any $n \ge d+1$. (Also observe that for any $X$ like above, $sa(X) \le a(X)$.)

In the second case, we consider the following example. Let us consider the cross-polytope $\mathrm{conv}\{\pm e_i\}$, where the $e_i$'s are the usual unit vectors. For $d+1 \le n \le 2d$ we write $n = d+k$, where $1 \le k \le d$. Then we take $\{\pm e_1, \ldots, \pm e_k\}$, whose convex hull is a $k$-dimensional cross-polytope. Then we add to it step by step $e_{k+1}, \ldots, e_d$. We have

$$sa(\{\pm e_1, \ldots, \pm e_k\}) = k$$

(taken in its affine hull $\mathbb{R}^k$). Addition of $e_{k+1}, \ldots, e_d$ means taking a $(d-k)$-fold pyramid over the above base. Applying Lemma 2.1, strictly antipodal case, $d-k$ times, we obtain

$$sa\left(\{\pm e_1, \ldots, \pm e_k, e_{k+1}, \ldots, e_d\}\right) =$$
$$k + 2k + (2k+1) + \ldots + (d+k-2) + (d+k-1) =$$
$$k + (d+k)(d+k-1)/2 - 2k(2k-1)/2.$$

Observe that both examples work for any $k \in \{1, \ldots, d\}$. □

**Proof of Theorem 2.4: 1.** We begin with the case $n \in [2d+1, 4d-1]$ odd. Let $P := \operatorname{conv} X$, that is a polytope with $|\operatorname{vert}(P-P)| = 2sa(X)$. Since $n = |\operatorname{vert} P|$ is odd, and for $Q := (P-P)/2$ we have that $|\operatorname{vert} Q|$ is even. Therefore $Q$ is not a translate of $P$. That is, the convex polytope $Q$ is reducible. Hence, by [36, Theorem 11], $2sa(X) = |\operatorname{vert} Q| \geq 4d$. Thus $sa(X) \geq 2d$.

**2.** For $\operatorname{conv} X$ a cross-polytope, evidently $sa(X) = d$. We turn to the case $n \in [d+2, 2d]$, $\operatorname{conv} X$ no cross-polytope. We use the same notations $P$ and $Q$ as above. Now we can have the following possibilities.
1) $|\operatorname{vert} P| \leq 2d-1$, while $|\operatorname{vert} Q| \geq 2d$. Then $Q$ is not a translate of $P$, so $Q$ is reducible, and like above, $sa(X) \geq 2d$.
2) $|\operatorname{vert} P| = 2d$, $P$ is no cross-polytope, while
A) $|\operatorname{vert} Q| \geq 2d+2$. Hence $Q$ is not a translate of $P$, so $Q$ is reducible, and $sa(X) \geq 2d$; or
B) $|\operatorname{vert} Q| = 2d < 4d$, hence $Q$ is a cross-polytope that is irreducible by [36, Theorem 11]. Hence $P$ is a translate of $Q$, thus $P$ is a cross-polytope, a contradiction.

**3.** It remains to give, for any $n \in [2d+1, 4d-1]$ odd, a set $X$ like in the theorem, with $sa(X) = 2d$. We begin by recalling the concept of *strongly isomorphic polytopes*, also called *analogous polytopes*. Two convex polytopes $P_1, P_2 \subset \mathbb{R}^d$ are strongly isomorphic if there is a combinatorial isomorphism between their face-lattices such that all corresponding face pairs are parallel and have the same orientation, cf. [29, p. 100] (there an equivalent definition is given). The positive linear combinations $\lambda_1 P_1 + \lambda_2 P_2$ of strongly isomorphic polytopes $P_1, P_2$ are strongly isomorphic to $P_1, P_2$. In particular, the sets of outward unit normals of the facets are the same. For $P_1, P_2, \lambda_1 P_1 + \lambda_2 P_2$ we have for the support functions

$$h_{\lambda_1 P_1 + \lambda_2 P_2}(u_i) = \lambda_1 h_{P_1}(u_i) + \lambda_2 h_{P_2}(u_i),$$

where the $u_i$'s are the outward unit normals of the facets. Moreover, this relation determines $\lambda_1 P_1 + \lambda_2 P_2$ uniquely among all positive linear combinations of $P_1, P_2$ (since the oriented facet hyperplanes determine a convex polytope).

Now we come to the example. Let us consider the regular cross-polytope $P_0 := \text{conv}\{\pm e_i\}$, where the $e_i$'s are the usual basic vectors. This can be considered as an antiprism with bases $\text{conv}\{e_i\}$ and $\text{conv}\{-e_i\}$. Since $d \geq 3$, we can omit the facet hyperplanes $\text{aff}\{e_i\}$ and $\text{aff}\{-e_i\}$, and then the remaining facet hyperplanes, with their original orientation, still will bound a convex polytope $P_0'$. Alternatively, we may describe this polytope by adding suitable simplices having one face $\text{conv}\{e_i\}$, or $\text{conv}\{-e_i\}$, respectively. (These simplices are regular pyramids with these bases.) These simplices can be described explicitly. We consider the facet $\text{conv}\{e_i\}$. All facets of $P_0$ are of the form $\text{conv}\{\varepsilon_i e_i\}$, where $\varepsilon_i \in \{-1, 1\}$. The facets adjacent to $\text{conv}\{e_i\}$ (i.e., having a common $(d-2)$-face with it) are therefore of the form $\text{conv}\{\varepsilon_i e_i\}$, where exactly one $\varepsilon_i$ is $-1$, the others are $1$. The corresponding facet hyperplanes have equations $\sum \varepsilon_i x_i = 1$, with $\varepsilon_i$ like above. Adding all $d$ such equations, we get $(d-2) \sum x_i = d$. Further, rewriting the equations as $-2x_i + x_1 + \ldots + x_d = 1$, we get $x_i = 1/(d-2)$. Clearly, the point $v := (1/(d-2), \ldots, 1/(d-2))$ satisfies also all other inequalities $\sum \varepsilon_i x_i \leq 1$, where at least two $\varepsilon_i$'s are $-1$. Hence, $\pm v \in \text{vert } P_0'$. Then $P_0' = \text{conv}\{\pm e_i, \pm v\}$, since, by symmetries, either all or no $\pm e_i$'s are vertices of $P_0'$ (and $P_0' = \text{conv}\{\pm v\}$ is impossible). Thus the added regular pyramids are $P_\pm := \pm \text{conv}\{e_i, v\}$.

Now let $H_\pm$ be parallel hyperplanes, properly intersecting each open lateral edge $\pm(e_i, v)$ of $P_\pm$. We will choose $H_\pm$ almost parallel to $\text{aff}\{e_i\}$. If they are exactly parallel, then, say $H_+$ has equation $\sum x_i = c_+$, where $c_+ > 1$. For $(x_1, \ldots, x_d) = -e_i$ we have $\sum x_i = -1$. So $H_+$ strictly separates $v$ from all $-e_i$, and this remains true if $H_+$ is only almost parallel to $\text{aff}\{e_i\}$. The analogous statement holds for $H_-$. Therefore, $H_+$ and $H_-$ cut off some simplices from $P_0'$, and actually, from $P_+$ and $P_-$. This truncation of $P_0'$ will be denoted by $P$.

Clearly, for any choice of the above hyperplanes $H_\pm$, with fixed outward normal unit vectors, all these polytopes $P$, as well as all the polytopes $-P$, are strongly isomorphic polytopes. In particular, also the polytopes $(P-P)/2$ are strongly isomorphic to the above ones.

Moreover, $(P - P)/2$ can be easily described among all these strongly isomorphic polytopes: it suffices to give the values of its support function at the outer unit normals of the facets. Observe that $P$ has $2^d + 2$ facet hyperplanes, namely all aff$\{\varepsilon_i e_i\}$, except $\pm$aff$\{e_i\}$, and $H_\pm$. All outer unit normals of these facet hyperplanes are fixed (recall that $H_\pm$ has a fixed outward normal unit vector). Also $-P$ has the same outer unit normals. Hence the strongly isomorphic polytope $(P-P)/2$ has the same outer unit normals. The support functions of $P$ and $-P$ coincide for the outward unit normals of the $2^d - 2$ facet hyperplanes aff$\{\varepsilon_i e_i\}$, hence also $(P - P)/2$ has equal support function values of the outward unit normals of these facet hyperplanes. In other words, these facet hyperplanes of $P$ and $(P - P)/2$ coincide. The support function of $(P - P)/2$ at the outer unit normal vector of the facet hyperplane $H_\pm$ is the arithmetic mean of the support functions of $P$ at the outer normal unit vectors of $H_+$ and $H_-$. More geometrically,

$$(P - P)/2 \text{ is bounded by all facet hyperplanes aff}\{\varepsilon_i e_i\},$$
$$\text{except } \pm \text{conv}\{e_i\}, \text{ and by the mid-hyperplane of the parallel}$$
$$\text{strip bounded by } H_+ \text{ and } -H_-, \text{ and by the mid-hyperplane}$$
$$\text{of the parallel strip bounded by } H_- \text{ and } -H_+. \quad (21)$$

We determine the vertices of $P$ and of $(P - P)/2$. At each vertex $\pm e_i$ one facet hyperplane was omitted. E.g., at $e_1$, the facet hyperplane aff$\{e_i\}$ was omitted. However, there remained still $2^{d-1} - 1$ facet hyperplanes at $e_1$, namely aff$\{e_1, \pm e_2, \ldots, \pm e_d\}$, where not all $\pm$ signs can be $+$ signs.

Their outer normal vectors are $(1, \pm 1, \ldots, \pm 1)$, where not all $\pm$ signs can be $+$ signs. We are going to show that these outer normal vectors linearly span $\mathbb{R}^d$; this will show that the considered vertex $e_1$ of $P_0$ will remain a vertex of $P$, as well.

Since $d \geq 3$, we have $(1, 1, 1, \ldots) + (1, -1, -1, \ldots) = (1, 1, -1, \ldots) + (1, -1, 1, \ldots)$, where the three dots denote 1's. Hence the linear hull of all $(1, \pm 1, \ldots, \pm 1)$, where not all $\pm$ signs can be $+$ signs, equals the linear hull of all $(1, \pm 1, \ldots, \pm 1)$, that further equals the linear hull of all $(\pm 1, \pm 1, \ldots, \pm 1)$, i.e., $\mathbb{R}^d$.

Therefore, all $\pm e_i$'s are vertices of $P$, and of $(P - P)/2$. Additionally, where $H_\pm$ (or the respective mid-hyperplane of a parallel strip) properly intersects the open segments $\pm(e_i, v)$, there are further vertices. Their total number is $4d$.

Now we consider the limiting situations, where $H_\pm$ may contain $\pm v$, or several $\pm e_i$'s. Recall that Minkowski sums are continuous in the summands, hence the above description (21) of $(P - P)/2$ will remain valid. We will avoid that both $H_+$ and $H_-$ should degenerate so that both contain $v$ and $-v$, or several $e_i$'s and several $-e_i$'s. Then the mid-hyperplane of the parallel strip $H_+, -H_-$ will strictly separate $v$ and the $e_i$'s. Therefore

$$|\text{vert}\,[(P-P)/2]| = 4d \text{ will continue to hold.}$$

At the same time, degeneration decreases $|\text{vert}\,P|$. If we choose a translation of $H^+$ (almost parallel to $\text{aff}\,\{e_i\}$) so that $H_+$ contains $v$, then $H_+$ contributes just one vertex to $P$, namely $v$. Then we may choose $H_-$ as strictly separating all $-e_i$ and $-v$. Then $H_-$ contributes $d$ vertices to $P$. However, we may prescribe that $H_-$ should contain exactly $k$ of the points $-e_i$, where $0 \le k \le d$ is arbitrary (and $H_\pm$ are almost parallel to $\text{aff}\,\{e_i\}$). Then $H_-$ contributes exactly $d - k \in [0, d]$ new vertices to $P$. Then

$$|\text{vert}\,P| = 2d + 1 + d - k \in [2d + 1, 3d + 1],$$

and here $|\text{vert}\,P|$ can be any number in $[2d+1, 3d+1]$.

In a similar way, let $H_+$ strictly separate all $e_i$'s and $v$, while $H_-$ contains exactly $k$ of the points $-e_i$, where $0 \le k \le d$ is arbitrary (and $H_\pm$ are almost parallel to $\text{aff}\,\{e_i\}$). Then, similarly as above,

$$|\text{vert}\,P| = 2d + d + d - k \in [3d, 4d],$$

and here $|\text{vert}\,P|$ can be any integer in $[3d, 4d]$.

In particular, we have an example with $sa(X) = 2d$, for all odd numbers in $[2d+1, 4d-1]$. □

**Proof of Proposition 2.5: 1.** Application of Lemma 2.1, strictly antipodal case, settles cases 1), 2), 3).

**2.** There remains the case of a triangular bipyramid $P$. This may have ten strictly antipodal vertex pairs, cf. [15]. Also, it may have seven strictly antipodal vertex pairs. We let the base be a regular triangle, and both pyramids be even (regular) pyramids. One is a regular tetrahedron, and the other has height some small $\varepsilon > 0$.

On the other hand, there are always at least seven strictly antipodal vertex pairs. Let $x_1, x_2$ be the two apices, and $x_3, x_4, x_5$ be the vertices of the base. Then $x_1, x_2$, and any pair among $x_3, x_4, x_5$ are strictly antipodal. In [21, p. 188] it was shown that, e.g., $x_1, x_3$ and $x_2, x_3$ are strictly antipodal, or, e.g., $x_1, x_3$ and $x_1, x_4$ are strictly antipodal. Thus we have already six strictly antipodal pairs, and we have to find still one.

Here we make some explanations to [21, p. 188].

**2A.** $P$ does not have collinear edges by strictly convex position, and does not have parallel non-collinear edges, since then some four vertices would be coplanar, while $P$ is a triangular bipyramid. The projection along a direction of general position in aff $\{x_3, x_4, x_5\}$ that we consider as horizontal different from those of the sides of conv $(x_3, x_4, x_5)$, can produce a quadrangle with a parallel pair of sides only if the following holds. The edges of $P$, being the inverse images of the parallel edges and therefore of the form $x_i x_j$, where $i \in \{1, 2\}$ and $j \in \{3, 4, 5\}$, have a normal transversal orthogonal to the projection direction. This excludes exactly one projection direction in the horizontal plane, unless the normal transversal of, e.g., aff $\{x_1, x_3\}$ and aff $\{x_2, x_4\}$ would be orthonormal to the horizontal plane, i.e., would be vertical. Then aff $\{x_1, x_3\}$ and aff $\{x_2, x_4\}$ would be horizontal, that is impossible.

Now let us suppose that, e.g., $x_1, x_3$ and $x_2, x_3$ are strictly antipodal, or that, e.g., $x_2, x_3$ and $x_1, x_4$ are strictly antipodal. In both cases, e.g., $x_1, x_3$ are strictly antipodal. Then choose a projection direction of general position in the horizontal plane, such that the projection of $x_3$ should lie strictly between the projections of $x_4$ and $x_5$. Then the image of $P$ by this projection is a quadrangle $Q$. By the general position of the projection direction, $Q$ has no parallel sides, and hence has two neighbouring sides whose endpoints are strictly antipodal vertex pairs of $Q$. By Lemma 2.2, the inverse images of these edges by the projection are strictly antipodal vertex pairs of $P$. Moreover, $x_3$ is not a vertex in any of these vertex pairs, since its projection is in the interior of $Q$. So its projection cannot form a strictly antipodal vertex pair with any point of $Q$. □

**Proof of Proposition 2.6:** The statement is a special case of the last statement of Theorem 2.4. □

**Proof of Theorem 2.7: 1.** Let $x \in X$, and let us project $x$ along a line in general position, close to the line connecting $x$ and the barycentre of $P := \operatorname{conv} X$. Let the image of $X$, or $P$, by this projection be $Y$, or $Q$, respectively. Then $|\operatorname{vert} Q| < |X| = d+k$, since the image of $x$ is in $\operatorname{rel int} Q$. Then for $d-1+k' := sa(Y) = sa(\operatorname{vert} Q)$ we have $1 \leq k' \leq k \leq d-1$.

By Lemma 2.2, we have $sa(X) \geq sa(Y)$. Taking into account that $x$ belongs to some strictly antipodal pair of $X$, which does not project to a strictly antipodal pair of $Y$, we have even

$$sa(X) \geq sa(Y) + 1,$$

where $1 \leq k' \leq k$. Taking minima for all $Y$, we obtain

$$sa(X) \geq \min\{(\underline{sa}^s_{d-1})'(d-1+k') \mid 1 \leq k' \leq k\} + 1, \qquad (22)$$

unless

$$k = d-1 \text{ and } Q \text{ is a } (d-1)\text{-dimensional cross-polytope}. \qquad (23)$$

We are going to exclude the case (23), unless we have (17) of the Theorem. In case (23) we would have, for a suitable numeration of the vertices of $Q$, that $y_1 + y_d = y_2 + y_{d+1} = \ldots = y_{d-1} + y_{2d-2}$. Let $x_i \in X \setminus \{x\}$ be the inverse image of $y_i$ by our projection. Then, for any $i, j$, we have that $x_i + x_{i+d-1} - x_j - x_{j+d-1}$ either equals 0, or is a non-zero vector, parallel to the line along which we projected. The second possibility can be avoided if the line, along which we projected, is not parallel to any non-zero vector $x_\alpha + x_\beta - x_\gamma - x_\delta$, for any different indices $\alpha, \beta, \gamma, \delta \in \{1, \ldots, d+k\}$. Thus there remains the case when all vectors $x_i + x_{i+d-1} - x_j - x_{j+d-1}$ are 0. Then the convex hull of all $x_i$'s being inverse images of the $(2d-2)$ $y_i$'s, are equal, and thus $Q$ is a $(d-1)$-dimensional cross-polytope. Still we have the point $x \in X$ used for choosing our projection. By $\operatorname{aff} X = \mathbb{R}^d$ this point $x$ does not belong to the affine hull of the former $(2d-2)$ $x_i$'s. Hence $P$ is a pyramid over a $(d-1)$-dimensional cross-polytope. By Lemma 2.1, $sa(X) = 3(d-1)$, i.e., (17) of the Theorem holds.

**2.** Suppose that we have

$$sa(X) \geq c \geq \min\{(\underline{sa}^s_{d-1})'(d-1+k') \mid 1 \leq k' \leq k\},$$

while we do not have (17) of the Theorem. Consider the graph with vertex set $X$, and edges the strictly antipodal pairs of $X$. This graph has average degree at least $2c/(d+k)$, hence for $x$ we can choose a vertex of degree at least $\lceil 2c/(d+k) \rceil$. Observe that strictly antipodal pairs of $X$, containing $x$, do not project to strictly antipodal pairs of $Y$. Then, like in **1**, we have, rather than (22),

$$sa(X) \geq sa(Y) + \lceil 2c/(d+k) \rceil \,.$$

Then

$$sa(X) \geq \min\{(\underline{sa}_{d-1}^s)'(d-1+k') \mid 1 \leq k' \leq k\} + \lceil 2c/(d+k) \rceil \,,$$

showing the first statement of the theorem.

**3.** Again suppose that we do not have (17) of the Theorem. For simplicity, we write

$$c_0 := \min\{(\underline{sa}_{d-1}^s)'(d-1+k') \mid 1 \leq k' \leq k\} \,.$$

Then (17) from **1** of the proof says

$$sa(X) \geq c_0 + 1 \,.$$

Then the last inequality from **2** implies

$$sa(X) \geq c_0 + 2(c_0 + 1)/(d+k) = c_0 \left[1 + 2/(d+k)\right] + 2/(d+k) \,.$$

That is, we have a better lower estimate for $sa(X)$ than $sa(X) \geq c_0 + 1$. Now, letting

$$c := c_0 \left[1 + 2/(d+k)\right] + 2/(d+k) \,,$$

the last inequality of **2** implies a further improvement, namely

$$sa(X) \geq c_0 + [2/(d+k)] \left[c_0 \left(1 + 2/(d+k)\right) + 2/(d+k)\right)$$
$$= c_0 \left[1 + 2/(d+k) + 2^2/(d+k)^2\right] + 2^2/(d+k)^2 \,.$$

Iteration of this process gives, by induction, lower estimates

$$sa(X) \geq c_0 \left[1 + 2/(d+k) + \ldots + 2^n/(d+k)^m\right] + 2^m/(d+k)^m,$$

and letting $m \to \infty$ we gain

$$sa(X) \geq c_0/[1 - 2/(d+k)],$$

that is equivalent to the last inequality of the Theorem. $\square$

**Proof of Proposition 2.8:** Let $X \in \mathcal{C}_d^s(d+k)$.

We apply Theorem 2.7. If we have (17) from Theorem 2.7, then $sa(X) = 3(d-1) \geq 2d + 1$, by $d \geq 4$. If we do not have (17) of Theorem 2.7, then

$$\begin{aligned}
sa(X) &\geq \lceil \min\left\{(\underline{sa}_{d-1}^s)'(d-1+k') \mid 1 \leq k' \leq k\right\} \cdot (d+k)/(d+k-2) \rceil \\
&\geq \lceil 2(d-1) \cdot (d+k)/(d+k-2) \rceil \geq \lceil 2(d-1) \cdot (2d-1)/(2d-3) \rceil \\
&= \lceil 2d + 2/(2d-3) \rceil = 2d + 1,
\end{aligned}$$

the second inequality following from Theorem 2.4, since $d - 1 \geq 3$. $\square$

**Proof of Theorem 2.9: 1.** The fact that for the given examples we have $sa(X) = 11$, follows from Theorem 2.3 and its proof. Observe that for $d = 4$ we have $(5d-2)/3 = 6$. Therefore the two examples in the proof of Theorem 2.3 give the same number of strictly antipodal pairs, for $n = 6$.

**2.** We turn to prove $\underline{sa}_4(6) \geq 11$.

Again we write $P := \text{conv} X$. First let $n = 6 = d + 2$. Then there is a Radon partition (i.e., there are two disjoint subsets of $X$ whose convex hulls intersect, i.e., $\sum_{i \in I} \lambda_i x_i = \sum_{j \in J} \mu_j x_j$ for some partition $\{I, J\}$ of $\{1, \ldots, 6\}$, and some $\lambda_i \geq 0$, $\mu_j \geq 0$ with $\sum_{i \in I} \lambda_i = \sum_{j \in J} \mu_j = 1$). Omitting all points from a Radon partition class for which the corresponding $\lambda_i$ or $\mu_j$ is 0, we obtain disjoint subsets $I, J \subset \{1, \ldots, 6\}$ such that $\sum_{i \in I} \lambda_i x_i = \sum_{j \in J} \mu_j x_j$, where $\lambda_i > 0$, $\mu_i > 0$, and $\sum_{i \in I} \lambda_i = \sum_{j \in J} \mu_j = 1$. Then, adding indices to $I \cup J$ means geometrically to take pyramids over the examples already constructed.

We have the following cases: The pair $(|I|, |J|)$ can be $(2, 2), (2, 3), (2, 4), (3, 4)$.

The case $(2,2)$ means a two-fold pyramid over a planar quadrangle. For a quadrangle the number of strictly antipodal vertex pairs is at least two, with equality only for a parallelogram. Now we apply Lemma 2.1, strictly antipodal case, twice. Therefore, in case $(2,2)$,

$$sa(X) \geq 2 + 4 + 5 = 11\,,$$

with equality if and only if the quadrangle is a parallelogram.

The case $(2,3)$ means a pyramid over a triangular bipyramid, whose vertex set is denoted by $Y$. By Lemma 2.2, $sa(Y) \geq 7$, hence $sa(X) \geq 7 + 5 = 12$.

There remain the non-pyramidal cases $(2,4)$ and $(3,3)$. Then the points in a partition class are vertices of simplices, of dimensions $1, 3$, or $2, 2$. Therefore each pair of points in the same partition class forms a strictly antipodal pair for the partition class, in its own affine hull. Hence they form strictly antipodal pairs of $X$, as well. Thus we have already

$$(2 \cdot 1)/2 + (4 \cdot 3)/2 = 7 \ \ \text{or} \ \ (3 \cdot 2)/2 + (3 \cdot 2)/2 = 6$$

strictly antipodal pairs in $X$. We still have to find strictly antipodal pairs $x_i, x_j$, with $i \in I$, $j \in J$.

We project $X$ along a line of general position in the affine hull of a not smaller partition class. The simplices $\operatorname{conv} \{x_i \mid i \in I\}$ and $\operatorname{conv} \{x_j \mid j \in J\}$ intersect in a single point, that is a relative interior point of both simplices. Their projections also lie in complementary subspaces and have a single point in common, that is a relative interior point of both of them. The dimensions of these subspaces are $1, 3 - 1$, or $2, 2 - 1$, respectively, i.e., in both cases $1, 2$. If we have dimensions $2, 2$, then the projection of the triangle, in whose affine hull we have the line along which we projected, will be a segment. In this case we suppose that the line is not parallel to any side of the respective triangle. Therefore

> some vertex of the triangle will project to a relative interior point of the projection of the triangle . (24)

If we have dimensions $1, 3$, then we may choose a line along which we project so that

> the projection of some vertex of the tetrahedron should be in the relative interior of the projection of the tetrahedron, that is a triangle in this case . (25)

Even for any sufficiently small perturbation of the line along which we project, in the affine hull of the respective partition class we will still have the same property (24) or (25), respectively. Therefore, the projection of $P$ is a triangular bipyramid, $Q$, say, whose vertex set we denote by $Y$. The set $Y$ is the projection of those five points of $X$ which are different from the vertices mentioned in (24) or (25), respectively.

By Proposition 2.5, we have $sa(Y) \geq 7$. Since the endpoints of an edge of the triangular base and of the diagonal are strictly antipodal in $Y$, these already give four strictly antipodal pairs in $Y$. We write $y_1, y_2$ for the two apices, and $y_3, y_4, y_5$ for the vertices of the triangular base. Therefore, there are at least $7 - 4 = 3$ strictly antipodal pairs $y_i, y_j$, where $i \in \{1, 2\}$ and $j \in \{3, 4, 5\}$. Using Lemma 2.2, the pairs of points of $X$, projecting to these pairs $y_i, y_j$, are strictly antipodal in $X$ and also lie in different partition classes $I, J$ of $X$. This implies three more strictly antipodal pairs in $X$, i.e., we have already $7 + 3 = 10$, or $6 + 3 = 9$, strictly antipodal pairs in $X$.

Let $x_i, x_j$ be one of these at least three strictly antipodal pairs in $X$. Then we may suppose that the vertex of a triangle mentioned in (24), or the vertex of the tetrahedron mentioned in (25), is in $\{x_i, x_j\}$, for a new projection. Then, using the new projection, we have 10, or 9, strictly antipodal pairs in $X$, but $x_i, x_j$ cannot be among them. Since one of these points projects to the relative interior of the respective projection (a segment, or a triangle), therefore we can add this pair $x_i, x_j$ to the ten, or nine, strictly antipodal pairs in $X$, found by this new projection, as a pair in different partition classes of $X$. That is, we have already $10 + 1 = 11$, or $9 + 1 = 10$, strictly antipodal vertex pairs.

In the first case, i.e., with cardinalities of partition classes $2, 4$, this gives the statement of the theorem. Thus there remains the second case, i.e., with cardinalities of partition classes $3, 3$.

Then we have found already six strictly antipodal pairs lying in the same partition class $I$ or $J$. Thus we have still $10 - 6 = 4$ strictly antipodal pairs in different partition classes. These pairs determine a two-partite graph $G$ with colour classes $I$ and $J$, and with the above said pairs as edges. Then the average degree in $G$ is $3 + 2/6 > 1$. Therefore a vertex with degree at least two must exist, say, $x_i$. (Observe that now $|I| = |J| = 3$, so it makes no difference whether $x_i \in I$ or $x_i \in J$.) Then repeat the second projection so that this vertex $x_i$ will play the role of the vertex in (24). (Observe that (25) is already excluded.) Then we have two strictly antipodal pairs $x_i, x_{j_1}$ and $x_i, x_{j_2}$, both of which we can add to the nine strictly antipodal vertex pairs. This gives the statement of the theorem for the case $|I| = |J| = 3$. This finishes the proof of the Theorem. □

**Proof of Proposition 2.10.** Let $X \in \mathcal{C}_d(d+k)$. Let in the graph $G$ with vertex set $X$ and edge set the strictly antipodal pairs in $X$, the maximal degree be $\ell$. Let this maximum be attained for $x_1 \in X$. Let us join $x_1$ with some interior point of $P$ by a line, and project $P$ along this line to the orthocomplementary linear $(d-1)$-subspace. Let this orthogonal projection map be $\pi$. Then $Y \in \mathcal{C}_{d-1}^s(|Y|)$. Also $Y \subset \pi X$, but $\pi x_1 \in \mathrm{relint}\, Q$, hence $d \le |Y| =: d + k' - 1 \le d + k - 1$, hence $1 \le k' \le k$.

Then by Lemma 2.2, each strictly antipodal pair in $Y$ lifts via $\pi$ to at least one strictly antipodal pair in $X$. However, the $\ell$ strictly antipodal pairs in $X$, containing $x_1$, are not lifted from any strictly antipodal pair in $Y$. Hence

$$sa(X) \ge sa(Y) + \ell \ge \min\{\underline{sa}_{d-1}^s(d-1-k') \mid 1 \le k' \le k\} + \ell,$$

as asserted. □

**Proof of Theorem 2.11.** Clearly it suffices to prove the theorem without the $\lceil$ and $\rceil$ signs.

By Proposition 2.10, and using that the maximal degree $\ell$ in $G$ is at least the average degree $2as(X)/(d+k)$, we get

$$sa(X) \ge \min\{\underline{sa}_{d-1}^s(d-1-k') \mid 1 \le k' \le k\} + \ell \ge$$
$$\min\{\underline{sa}_{d-1}^s(d-1-k') \mid 1 \le k' \le k\} + 2as(X)/(d+k).$$

Leaving out the middle term, we get

$$sa(X) \cdot (1 - 2/(d+k)) \ge \min\{\underline{sa}_{d-1}^s(d-1+k') \mid 1 \le k' \le k\},$$

which is equivalent to the first statement of Theorem 2.11.

We deduce the second statement of the theorem from the first one. Applying the first statement for $k' \in \{1,\ldots,k\}$, and then taking minimum for all these values $k'$, we get

$$\min\{\underline{sa}_d^s(d+k') \mid 1 \le k' \le k\} \ge$$
$$\min\{\bigl(\min\{\underline{sa}_{d-1}^s(d-1+k'') \mid 1 \le k'' \le k'\}\bigr) \mid 1 \le k' \le k\}$$
$$\cdot (d+k)/(d+k-2) = \min\{\underline{sa}_{d-1}^s(d-1+k'') \mid 1 \le k'' \le k\},$$

which implies the second statement of the theorem.

Conversely, supposing the second statement of the theorem, we have

$$\underline{sa}^s_d(d+k) \geq \min\{\underline{sa}^s_d(d+k') \mid 1 \leq k' \leq k\} \geq \\ \min\{\underline{sa}^s_{d-1}(d-1+k') \mid 1 \leq k' \leq k\} \cdot (d+k)/(d+k-2),$$

which implies the first statement of the theorem. $\square$

**Proof of Lemma 2.12.** Observe that $a$ is a common point of $P$ and $P + a - b$; suppose $a = 0$. Moreover, this is their only common point. Namely, both $P$ and $P - b$ are the unions of all their respective chords of the given direction. Two such respective chords can have a common point only if they are collinear. And if they are collinear, but different from $[a, b]$, then one of them is a translate of the other one by the vector $b$, while their length is by hypothesis less than $\|b\|$, so they are disjoint.

Then $0$, being the unique common point of $P$ and $P - b$, satisfies that the support cones $S_1$ and $S_2$ of $P$ and $P - b$ at $0$ have as intersection $\{0\}$. Then $S_1 \cap S^{d-1}$ and $S_2 \cap S^{d-1}$ are disjoint spherically convex spherical $(d-1)$-polytopes. Let their spherical distance be $\delta > 0$. The dual cones $S_i^*$ (normal cones of $P$ and $P - b$ at $0$) are $d$-dimensional. Choosing some inner points $u_i^* \in S^{d-1}$ of them, we see that $S_i \setminus \{0\}$ lies in an open halfspace bounded by the orthocomplement of $u_i^*$. Therefore $S_i \cap S^{d-1}$ lies in an open hemisphere of centre $-u_i^*$. Then central projection with centre $0$ establishes a bijective map $p$ from this open hemisphere to the tangent hyperplane of $S^{d-1}$ at $-u_i^*$. Both $p$ and $p^{-1}$ preserve convexity.

Consider, for some sufficiently small $\varepsilon > 0$, the $\varepsilon$-neighbourhood $[p(S_i \cap S^{d-1})]_\varepsilon$ of $p(S_i \cap S^{d-1})$ in this hyperplane, which is compact and convex. Then $p^{-1}\left[[p(S_i \cap S^{d-1})]_\varepsilon\right]$ is spherically convex, compact, has diameter less than $\pi$, and we may suppose that it lies in the spherical $\delta/3$-neighbourhood of $S_i \cap S^{d-1}$. Let the cone $S_i(\varepsilon)$, with apex $0$, be spanned by $p^{-1}\left[[p(S_i \cap S^{d-1})]_\varepsilon\right]$. Then $S_1(\varepsilon) \cap S_2(\varepsilon) = \{0\}$, hence $S_1(\varepsilon)$ and $S_2(\varepsilon)$ can be separated by a hyperplane $H$, containing $0$. Since $S_i \setminus \{0\} \subset \text{int } S_i(\varepsilon)$, therefore $H \cap S_i = \{0\}$. Hence $P \cap H = \{0\}$. Moreover, also $P \cap (H + b)$, which is a translate of $(P - b) \cap H$, consists of the single point $b$. That is, $a = 0$ and $b$ are strictly antipodal in $P$. $\square$

**Proof of Proposition 2.13. 1.** We write $n := d + k$ and $X = \{x_1, \ldots, x_n\}$. Let $x_1 \in X$ be a vertex of $G$ of maximal degree $\ell \leq n - 2$. Then there are, say, $x_2, \ldots, x_{\ell+1} \in X$ such that $\{x_1, x_2\}, \ldots, \{x_1, x_{\ell+1}\}$ are edges of $G$, but $\{x_1, x_n\}$ is not an edge of $G$.

For subsets $X', X'' \subset X$ we write $sa(X; X', X'')$ for the number of strictly antipodal pairs $\{x', x''\}$ of $X$, with $x' \in X'$ and $x'' \in X''$. We write $sa(X; X') := sa(X; X', X)$. For $X' = \{x'\}$ we write $sa(X; x') := sa(X; \{x'\})$.

Then $x_1, x_n \in X$ are not strictly antipodal, hence by Lemma 2.12, $[x_1, x_n]$ is not a unique maximal segment of its direction, contained in $P$. Hence there is a segment $[x'_1, x'_n] \subset P$, parallel to but not collinear with $[x_1, x_n]$, which has at least the same length as $[x_1, x_n]$. We may suppose that these lengths are equal, and that $x_1 x_n x'_n x'_1 \subset P$ is a parallelogram with this order of vertices. Then the half-lines $x_1 x'_1$ and $x_n x'_n$ leave $P$ at points $x_1^*$ and $x_n^*$ (different from $x_1$ and $x_n$).

Denote by $\pi$ the orthogonal projection of $\mathbb{R}^d$ along the line $x_1 x'_1$ (which we consider as vertical), to the orthocomplementary linear $(d-1)$-subspace. We write $m := |Y| \geq d$; then $Y \in \mathcal{C}^s_{d-1}(m)$.

For $y \in Y$, we have that $(\pi^{-1}y) \cap P$ is either a vertex, or an edge of $P$. (26)

By Lemma 2.2, if $\{y', y''\}$ is a strictly antipodal pair of $Y$, then $sa(Y; \pi^{-1}y', \pi^{-1}y'') \geq 1$, and if at least one of $(\pi^{-1}y') \cap P$ and $(\pi^{-1}y'') \cap P$ is an edge, then $sa(Y; \pi^{-1}y', \pi^{-1}y'') \geq 2$. If both of them are edges, then the two strictly antipodal pairs lifting $\{y', y''\}$ are disjoint (hence there are no more strictly antipodal pairs lifting $\{y', y''\}$).

**2.** We distinguish four cases.
(a) $\pi x_1, \pi x_n \notin Y$;
(b) $y_1 := \pi x_1, y_m =: \pi x_n \in Y$;
(c) $y_1 := \pi x_1 \in Y$ but $\pi x_n \notin Y$;
(d) $\pi x_1 \notin Y$ but $y_m =: \pi x_n \in Y$.

**3.** First we prove (1), i.e., that $m = |Y| \leq |X| - 2 = n - 2$.

In Case (a) (1) holds, since $P = \text{conv}\, X$ implies $Q = \text{conv}(\pi X) = \text{conv}\,(\pi(X \setminus \{x_1, x_n\}))$, hence $|Y| \leq |X| - 2$.

In Case (b), $(\pi^{-1}y_1) \cap P$ contains the non-trivial segment $[x_1, x'_1]$. Hence it cannot be a vertex, therefore by (26) it is an edge, namely $[x_1, x_1^*]$, of $P$. Similarly, $(\pi^{-1}y_m) \cap P$ is an edge, namely $[x_n, x_n^*]$ of $P$. Therefore in Case (b) (1) holds, since $P = \text{conv}\, X$ implies $Q = \text{conv}(\pi X) = \text{conv}\,(\pi(X \setminus \{x_1^*, x_n^*\}))$, hence $|Y| \leq |X| - 2$.

In Case (c), similarly as in Case (b), $(\pi^{-1} y_1) \cap P$ is an edge, namely $[x_1, x_1^*]$ of $P$. Therefore in Case (c) (1) holds, since $P = \operatorname{conv} X$ implies $Q = \operatorname{conv}(\pi X) = \operatorname{conv}(\pi(X \setminus \{x_1^*, x_n\}))$, hence $|Y| \leq |X| - 2$.

In Case (d), similarly as in Case (b), $(\pi^{-1} y_m) \cap P$ is an edge, namely $[x_n, x_n^*]$ of $P$. Therefore in Case (d) (1) holds, since $P = \operatorname{conv} X$ implies $Q = \operatorname{conv}(\pi X) = \operatorname{conv}(\pi(X \setminus \{x_1, x_n^*\}))$, hence $|Y| \leq |X| - 2$.

Thus (1) is proved for all of the cases (a), (b), (c), (d). Hence $Y \in \mathcal{C}_{d-1}^s(|Y|)$ and $d \leq |Y| \leq |X| - 2 \leq 2d - 3$. These imply the second inequality in (2).

**4.** Now we prove the first inequality in (2), namely that $sa(X) \geq sa(Y) + \ell$ (the second inequality in (2) was just shown, so this will prove (2)).

By Lemma 2.2 each strictly antipodal pair of $Y$ lifts by $\pi$ to at least one strictly antipodal pair of $X$. Therefore for the first inequality in (2) we have to find still $\ell$ additional strictly antipodal pairs of $X$, as compared to $Y$. We mean by this that either they are not lifted by $\pi$ from any strictly antipodal pair from $Y$, or are lifted doubly by $\pi$ from some strictly antipodal pair from $Y$. (Later "by $\pi$" will be omitted, but implied.)

**4a.** In Case (a) observe that each strictly antipodal pair from $Y$ lifts, by Lemma 2.2, to at least one strictly antipodal pair of $X$. By hypothesis, $x_1$ is contained in $sa(X; x_1) = \ell$ additional strictly antipodal pairs of $X$, namely in $\{x_1, x_2\}, \ldots, \{x_1, x_{\ell+1}\}$. Hence

$$sa(X) \geq sa(Y) + sa(X; x_1) = sa(Y) + \ell.$$

**4b.** We turn to case (b).

For $2 \leq i \leq \ell + 1$, if $\{y_1, \pi x_i\}$ are not strictly antipodal in $Y$ (possibly because they are equal), then the pairs $\{x_1, x_i\}$ are additional strictly antipodal pairs in $X$. The set of these $i$'s is denoted by $I_1$. *Thus we have obtained $|I_1|$ many additional strictly antipodal pairs of $X$.*

It remains to consider those $i \in \{2, \ldots, \ell+1\}$, for which $\{y_1, \pi x_i\}$ is strictly antipodal in $Y$ – thus $y_1 \neq \pi x_i$ and $y_1, \pi x_i \in Y$. The set of these $i$'s is denoted by $I_2$.

First we show that for different $x_i, x_j \in I_2$ also the points $\pi x_i, \pi x_j$ are different. In fact, if $\pi x_i = \pi x_j$ held, then $x_i, x_j$ would be different vertices of $P$ in $(\pi^{-1}(\pi x_i)) \cap P$, hence $(\pi^{-1}(\pi x_i)) \cap P$ would be an edge, namely $[x_i, x_j]$ of $P$ (recall $\pi x_i \in Y$ and (26)). Then $[x_1, x_1^*]$ and $[x_i, x_j]$ would be parallel but distinct edges of $P$ (recall $\pi x_1 = y_1 \neq \pi x_i$). By the second statement of Lemma 2.2, $\{y_1, \pi x_i = \pi x_j\}$ would lift to exactly two strictly antipodal pairs in $X$, which are moreover disjoint. However, $\{x_1, x_i\} \cap \{x_1, x_j\} \neq \emptyset$, a contradiction.

Therefore any strictly antipodal pair of $X$ lifting $\{y_1, \pi x_i\}$, and any strictly antipodal pair of $X$ lifting $\{y_1, \pi x_j\}$, for different $x_i, x_j \in I_2$, are different. Hence it suffices to show that any strictly antipodal pair $\{y_1, \pi x_i\}$ of $Y$, for $i \in I_2$, lifts to an additional strictly antipodal pair of $X$. Thus we will have $|I_2|$ many additional strictly antipodal pairs of $X$.

Since $(\pi^{-1}(y_1)) \cap P = [x_1, x_1^*]$ is an edge of $P$, therefore $\{y_1, \pi x_i\}$ lifts to two different strictly antipodal pairs of $X$, one of which is $\{x_1, x_i\}$. Hence here we have one additional strictly antipodal pair in $X$, as compared to $Y$, lifting $\{y_1, \pi x_i\}$ (and containing $x_1^*$). *Thus we have obtained $|I_2|$ many additional strictly antipodal pairs of $X$.*

Hence

$$sa(X) \geq sa(Y) + |I_1| + |I_2| = sa(Y) + sa(X; x_1) = sa(Y) + \ell.$$

**4c.** We turn to Case (c). Similarly as in Case (b), there are $|I_1| + |I_2| = sa(X; x_1) = \ell$ additional strictly antipodal pairs in $X$ (each containing either $x_1$ or $x_1^*$). Hence

$$sa(X) \geq sa(Y) + sa(X; x_1) = sa(Y) + \ell.$$

**4d.** We turn to case (d). Similarly as in Case (a), there are $sa(X; x_1) = \ell$ additional strictly antipodal pairs in $X$, namely $\{x_1, x_2\}, \ldots, \{x_1, x_{\ell+1}\}$. Hence

$$sa(X) \geq sa(Y) + sa(X; x_1) = sa(Y) + \ell.$$

Thus (2) is proved for all of the cases (a), (b), (c), (d). Since (1) has been already proved in **3**, the proof is finished. □

**Proof of Proposition 2.14.** The proof is the same as that of the first statement of Theorem 2.11, except that in place of Proposition 2.10 we use Proposition 2.13. □

**Proof of Theorem 2.15. 1.** We write $n := d + k$, and $X =: \{x_1, \ldots, x_n\}$. For $X$ the vertex set of a pyramid over a $(d-1)$-dimensional cross-polytope we have by Lemma 2.1 that $sa(X) = 3(d-1)$, hence $\underline{sa}_d^s(2d-1) \leq 3(d-1)$. Hence the inequality in Theorem 2.15 implies the equality in Theorem 2.15.

**2.** We turn to prove the inequality of Theorem 2.15. Observe that for $d = 2$ it is evident, while for $d = 3$ it follows from (9). We use induction for $d$: we suppose that the Theorem is valid for $\mathbb{R}^{d-1}$, where we will suppose that $d \geq 4$.

For $n = d+1$ we have $\underline{sa}_d^s(d+1) = d(d+1)/2 \geq 3(d-1)$ for $d \geq 4$. So we will suppose $n \geq d+2$.

Let $x_1 \in X$ be a vertex of $G$ of maximal degree; let this degree be $\ell = sa(X; x_1) \leq n - 1$. Then, say, $\{x_1, x_2\}, \ldots, \{x_1, x_{\ell+1}\}$ are edges of $G$, and $\{x_1, x_{\ell+2}\}, \ldots, \{x_1, x_n\}$ are not edges of $G$.

Then $\ell$ is at least the average degree of a vertex of $G$, therefore by Theorem 2.4, second statement

$$\ell \geq \lceil 2 \cdot 2d/n \rceil \geq \lceil 2 \cdot 2d/(2d-1) \rceil \geq 3. \tag{27}$$

**3.** Now we distinguish two cases.
(1) For each $x \in X \setminus \{x_1\}$ we have that $x_1, x$ are strictly antipodal (i.e., $\ell = n - 1$).
(2) We have that $x_1, x_n$ are not strictly antipodal (i.e., $\ell \leq n - 2$).

**4.** We begin with Case (1). By Proposition 2.10, with the notations from it, we have

$$Y \in \mathcal{C}_{d-1}^s(|Y|), \ |Y| \leq |X| - 1 \leq 2d - 2, \text{ and } sa(X) \geq sa(Y) + \ell.$$

If here we have
(1a) $|Y| \leq 2d - 3 = 2(d-1) - 1$,
then with the induction hypothesis, and by $n \geq d + 2 \geq 4$, we get

$$sa(X) \geq sa(Y) + \ell = sa(Y) + n - 1 \geq 3(d-2) + (n-1) \geq 3(d-2) + 3 = 3(d-1).$$

There remains the case
(1b) $|Y| = 2d - 2$.
Then from above also $\ell = n - 1 = 2d - 2$ follows, and hence by Proposition 2.10 and (7) we have

$$sa(X) \geq sa(Y) + \ell \geq \lceil (2d-2)/2 \rceil + (2d-2) = 3(d-1).$$

This finishes the proof in Case (1).

**5.** We turn to the proof in Case (2). By Proposition 2.13, the induction hypothesis and (27) we get

$$sa(X) \geq \min\{\underline{sa}^s_{d-1}(d-1+k') \mid 1 \leq k' \leq k-1\} + \ell \geq 3(d-2) + 3 = 3(d-1).$$

This finishes the proof in Case (2), and hence that of the theorem. □

**Proof of Theorem 2.16.** Clearly it suffices to prove the theorem without the ⌈ and ⌉ signs.

We denote the graph with vertex set $X$ and edge set the strictly antipodal pairs of $X$ by $G$. The maximal degree in $G$ is denoted by $\ell \leq d + k - 1$. We have the same distinction of Cases (1) and (2) as in **3** of the proof of Theorem 2.15.

In Case (1) we apply Proposition 2.10, (2), with $\ell = d + k - 1$. Then we obtain as a lower bound the expression in the second line of (19). In Case (2) we apply Proposition 2.14, whose hypothesis $\ell \leq d + k - 2$ is satisfied by Case (2). Then we obtain as a lower bound the expression in the third line of (19).

In both of these cases a valid lower bound is the minimum of the above two lower bounds, which proves the theorem. □

**Proof of Lemma 2.17.** We have to prove

$$f(d-1, k) + (d+k-1) \geq f(d-1, k-1) \cdot (1 + 2/(d+k-2)),$$

i.e.,

$$f(d-1, k) - f(d-1, k-1) + (d+k-1)$$
$$\geq f(d-1, k-1) \cdot 2/(d+k-2).$$

Here the left hand side equals $d$, and the right hand side equals $d - k + 1$, so the last inequality holds, which proves the lemma. □

**Proof of Theorem 2.18.** The $\Theta(\cdot)$ relation in the theorem follows from $(d-k)\Theta(d) + O(d) = \Theta((d-k)d)$.

We use the notation $f(d, k)$ from Lemma 2.17. Thus we have to prove $\underline{sa}^s_d(d+k) \geq f(d, k)$.

For $k = 1$ we have $\underline{sa}_d^s(d+1) = d(d+1)/2 = f(d,1)$. For $k = d-1$ we have by Theorem 2.15 and $d \geq 2$ that $\underline{sa}_d^s(2d-1) = 3(d-1) \geq 2d-1 = f(d, d-1)$. Observe that these already prove all cases $d \leq 3$ and $1 \leq k \leq d-1$.

Now we make induction for $d$. We suppose $d \geq 4$, and $2 \leq k \leq d-2$. Then we may apply Theorem 2.16. Since $f(d,k)$ is a decreasing function of $k \in \{1, \ldots, d-1\}$, (19) will take the simpler form

$$\underline{sa}_d^s(d+k) \geq \min\{\underline{sa}_{d-1}^s(d-1+k) + (d+k-1),$$
$$\underline{sa}_{d-1}^s(d-1+k-1) \cdot (d+k)/(d+k-2)\}.$$

By the induction hypothesis, the right hand side here can be further estimated below by

$$\min\{f(d-1, k) + (d+k-1), f(d-1+k-1) \cdot (d+k)/(d+k-2)\}.$$

However, by Lemma 2.17 (which applies by $2 \leq k \leq d-2$), the minimum in the last formula equals

$$f(d-1+k-1) \cdot (d+k)/(d+k-2) =$$
$$[(d-1+k-1)(d-k+1)/2] \cdot (d+k)/(d+k-2) = f(d,k),$$

proving the theorem. □

**Proof of Corollary 2.19.** We know $\underline{sa}_5^s(6) = 15$, and Theorem 2.15 implies $\underline{sa}_5^s(9) = 12$. Recall that we know $\underline{sa}_4^s(4+k)$ for each $1 \leq k \leq 3$, cf. Theorem 2.9 and Theorem 2.15. Then (19) gives $\underline{sa}_5^s(7) \geq 14$ and $\underline{sa}_5^s(8) \geq 14 + 2/3$, and (18) gives $\underline{sa}_5^s(7) \geq 15.4$ and $\underline{sa}_5^s(8) \geq 12$. Also using the upper bounds from Theorem 2.3, we get $\underline{sa}_5^s(7) = 16$ and $15 \leq \underline{sa}_5^s(8) \leq 17$. □

**Proof of Proposition 2.20.** The $\Theta(\cdot)$ relation follows from Theorems 2.3 and 2.18.

For $2 \leq k \leq (2d-2)/3$ the formula in the upper bound of $\underline{sa}_d^s(d+k)$ increases, and the formula in its lower bound decreases, hence their quotient is maximized for $k = (2d-2)/3$ (although this may not be an integer). Taking in account that for $k = (2d-2)/3$ the formulas in examples (1) and (2) in Theorem 2.3 give the same value, it suffices to investigate the case $(2d-2)/3 \leq k \leq d-1$. Then we write $\ell := d - k$. Then we have to show that the upper bound is at most 2 times the lower bound. This is a quadratic inequality for $\ell$, and when written out and rearranged, it becomes

$$\ell^2 + \ell + 2d \geq 0,$$

which holds.

As stated in the Remarks after the respective theorems, for $x \in [2/3, 1)$ and $k = \lfloor xd \rfloor$, we have that for $(2d-2)/3 \leq k \leq d-1$ the upper bound is $(d^2/2) \cdot 3(1-x)(x+1/3) + O(d)$, and the lower bound is $(d^2/2)(1-x)(1+x) + O(d)$. For $x \nearrow 1$ their quotient is $(1+3x)/(1+x) + O(1/d)$, and here $(1+3x)/(1+x) \nearrow 2$. $\square$

**Proof of Theorem 2.21.** The statement follows from Theorems 2.15 and 2.16. $\square$

# 5 A large strictly antipodal set

**Proof of Theorem 3:** We mimic the proof of [2, Lemma 9.11.2], with a small modification that will increase the base of the exponential given there.

By Csikós-Kiss-Swanepoel-de Wet [8] there exist four sets $S_1, S_2, S_3, S_4 \in \mathbb{R}^3$ that are weakly strictly antipodal (i.e., for any distinct $i, j \in \{1, \ldots, 4\}$, and any $x_i \in S_i, x_j \in S_j$ there are distinct parallel supporting planes of $\text{conv}\,(\bigcup_{i=1}^{4} S_i)$, one containing $x_i$, the other containing $x_j$, and $(\bigcup_{i=1}^{4} S_i) \setminus \{x_i, x_j\}$ lying strictly between these planes - that implies that $S_1, S_2, S_3, S_4$ are pairwise disjoint), where $S_1, S_2, S_3$ are embedded $C^1$ arcs, and $S_4$ is a singleton. (We remark that a possibly simpler example is obtained from the construction in [17, p. 464]. We take $S_1, S_2, S_3$ as circular arcs of small curvature, joining the endpoints of three mutually skew edges of $[-1, 1]^3$, each arc lying in the 2-plane spanned by the respective edge and 0, and turning with its concave side to 0. Then $S_4$ can be chosen as any of the two vertices of $[-1, 1]^3$, not contained in any of the three skew edges.)

Now let $k := \lfloor d/3 \rfloor$. We consider $\mathbb{R}^{3k}$ as embedded in $\mathbb{R}^d$ as the $x_1 \ldots x_{3k}$ coordinate subspace. Our constructed set $X$ will be in $\mathbb{R}^{3k}$. For $d \equiv 1 \pmod{3}$ we take the vertices of a pyramid over $\text{conv}\,X$, and for $d \equiv 2 \pmod{3}$ we take the vertices of a two-fold pyramid over $\text{conv}\,X$. By Lemma 2.1, strict antipodality of $X$ implies strict antipodality of the respective vertex set, and its cardinality is $|X| + 1$, or $|X| + 2$, respectively.

From now on we work in $\mathbb{R}^{3k}$, that we write as $\mathbb{R}^3 \oplus \ldots \oplus \mathbb{R}^3$ ($k$ direct summands). We write $(\mathbb{R}^3)_l$ for the $l$'th direct summand, and $\pi_l$ for the orthogonal projection of $\mathbb{R}^{3k}$ to $(\mathbb{R}^3)_l$, for $1 \leq l \leq k$. In the $l$'th summand, we define $\sigma_{l1}, \sigma_{l2}, \sigma_{l3}, \sigma_{l4}$ as the images of $\sigma_1, \sigma_2, \sigma_3, \sigma_4$ by the linear map of $\mathbb{R}^3$ to $(\mathbb{R}^3)_l$, mapping the usual basic unit vectors $\mathbf{e}_1, \mathbf{e}_2, \mathbf{e}_3$ of $\mathbb{R}^3$ to the basic unit vectors $\mathbf{e}_{3l-2}, \mathbf{e}_{3l-1}, \mathbf{e}_{3l}$ of $(\mathbb{R}^3)_l \subset \mathbb{R}^{3k}$.

Now let $(i_1, \ldots, i_k) \in \{1, 2, 3, 4\}^k$. We choose points $x_{i_1 \ldots i_k} \in \sigma_{1 i_1} \oplus \ldots \oplus \sigma_{k i_k}$, altogether $4^k$ points, in such a way that for $1 \leq l \leq k$ and $i_l \in \{1, 2, 3\}$, all $\pi_l(x_{i_1 \ldots i_l \ldots i_k})$, that lie in $\sigma_{l i_l}$, are different. For $i_l = 4$, of course, all $\pi_l(x_{i_1 \ldots i_l \ldots i_k})$, that lie in $\sigma_{l4}$, coincide.

We claim that $X := \{x_{i_1 \ldots i_k} \mid (i_1, \ldots, i_k) \in \{1, 2, 3, 4\}^k\}$ is a strictly antipodal set. Observe that for $(i_1, \ldots, i_k) \neq (j_1, \ldots, j_k)$ we have $x_{i_1 \ldots i_k} \neq x_{j_1 \ldots j_k}$, since if $i_l \neq j_l$, then $\pi_l(x_{i_1 \ldots i_k}) \in \sigma_{l i_l}$ and $\pi_l(x_{j_1 \ldots j_k}) \in \sigma_{l j_l}$, and $\sigma_{l i_l}$ and $\sigma_{l j_l}$ are disjoint. Therefore $|X| = 4^k$.

Now let $(i_1, \ldots, i_k) \neq (j_1 \ldots j_k)$. We claim that $x_{i_1 \ldots i_k}$ and $x_{j_1 \ldots j_k}$ are strictly antipodal in $X$. We may suppose $i_1 \neq j_1$. Moreover, for $i_1, j_1 \in \{1, 2, 3\}$ we may suppose $i_1 = 1, j_1 = 2$, and for $\{i_1, j_1\} \not\subseteq [1, 3]$ we may suppose $i_1 = 1, j_1 = 4$.

Then $\pi_1(x_{i_1 \ldots i_k}) \in \sigma_{1 i_1}$ and $\pi_1(x_{j_1 \ldots j_k}) \in \sigma_{1 j_1}$ are strictly antipodal with respect to $\sigma_{11} \cup \sigma_{12} \cup \sigma_{13} \cup \sigma_{14}$, hence also with respect to $\pi_1(X) \subset \sigma_{11} \cup \sigma_{12} \cup \sigma_{13} \cup \sigma_{14}$. Then there exists a parallel slab $S \subset (\mathbb{R}^3)_1$, such that $\pi_1(x_{i_1 \ldots i_k})$ and $\pi_1(x_{j_1 \ldots j_k})$ belong to different boundary planes of $S$, and all other points of $\pi_1(X)$ lie in rel int $S$.

Now consider $\pi_1^{-1}(S)$. This is a parallel slab in $\mathbb{R}^{3k}$, and by $\pi_4(X) \subset S$ we have $X \subset \pi_1^{-1}(S)$. Moreover, $x_{i_1 \ldots i_k}$ and $x_{j_1 \ldots j_k}$ belong to different boundary hyperplanes of $\pi^{-1}(S)$ and, for any $x_{\alpha_1 \ldots \alpha_k} \in X \setminus \{x_{i_1 \ldots i_k}, x_{j_1 \ldots j_k}\}$ we have either $x_{\alpha_1 \ldots \alpha_k} \in \text{int } S$, or $x_{\alpha_1 \ldots \alpha_k} \in \pi_1^{-1}\{\pi_1(x_{i_1 \ldots i_k}), \pi_1(x_{j_1 \ldots j_k})\}$, i.e.,

$$\pi_1(x_{\alpha_1 \ldots \alpha_k}) \in \{\pi_1(x_{i_1 \ldots i_k}), \pi_1(x_{j_1 \ldots j_k})\}. \qquad (28)$$

Now suppose $i_1 = 1, ..., j_1 = 2$. Then $\pi_1(x_{i_1 \ldots i_k}) \in \sigma_{1 i_1} = \sigma_{11}$, and $\pi_1(x_{j_1 \ldots j_k}) \in \sigma_{1 j_1} = \sigma_{12}$. However, this is impossible, since $X$ was chosen so that whenever $\pi_l(x_{\alpha_1 \ldots \alpha_k}) \in \sigma_{l i_l}$, with $i_l \in [1, 3]$, then all these $\pi_l(x_{\alpha_1 \ldots \alpha_l})$'s should be different.

Next suppose $i_1 = 1$ and $j_1 = 4$. Then $\pi_1(x_{i_1...i_k}) \in \sigma_{1i_1} = \sigma_{11}$, and $\pi_1(x_{j_1...j_k}) \in \sigma_{1j_2} = \sigma_{14}$. Like above, $\pi_1(x_{\alpha_1}...x_{\alpha k}) = \pi_1(x_{i_1}...x_{i_k}) \in \sigma_{1i_1} = \sigma_{11}$ is impossible. There remains the case that $\pi_1(x_{\alpha_1...\alpha_k}) = \pi_1(x_{j_1...j_k}) \in \sigma_{1j_1} = \sigma_{14}$. Since $|\sigma_{14}| = 1$, this surely happens when $\alpha_1 = 4$, and since $\sigma_{11}, \sigma_{12}, \sigma_{13}, \sigma_{14}$ are pairwise distinct, this does not happen for $\alpha_1 \in [1, 3]$. That is, there remains exactly the case when $\alpha_1 = 4$.

Let $Y := \{x_{\alpha_1...\alpha_k} \mid \alpha_1 = 4\} \ni x_{j_1...j_k}$. If we knew that $Y$ is in strictly convex position, in aff $Y = \sigma_{14} \oplus (\mathbb{R}^3)_2 \oplus ... \oplus (\mathbb{R}^3)_4$, then we could perturb $S$ a bit, so that it still contains $X$, one of its boundary hyperplane $H_1$ contains the only point $x_{i_1...i_k}$ of $X$, and the other boundary hyperplane $H_2$ rotates about a $(d-2)$-plane contained in it, that exposes the point $x_{j_1...j_k} = x_{4j_2...j_k}$ of the set $\{x_{\alpha_1...\alpha_k} = x_{4\alpha_2...\alpha_k}\}$, in the original position of $H_2$, and such that int $S$ contains all $x_{\alpha_1...j_k} = x_{4\alpha_2...\alpha_k}$ different from $x_{j_1...j_k}$. This would show that $X$ is strictly antipodal.

However, $Y$ is nothing else than our example $X$, with $k$ replaced by $k-1$ (multiplied by the singleton $\sigma_{14}$). This makes it possible to use induction for $k$. For $k = 1$, $X$ is strictly antipodal. If for $k-1$ the set $Y$ is strictly antipodal, then it is in strictly convex position, and the argument above shows that $X$ is strictly antipodal as well. $\square$

**Remark.** The result announced by [6, p. 271] was obtained by a 4-dimensional construction, and taking the $\lfloor d/4 \rfloor$'th power of it. This construction was the following. We take three mutually skew edges of $[-1, 1]^3$. Then, embedding $[-1, 1]^3$ to $[-1, 1]^4$ as $[-1, 1]^3 \times \{-1\}$, we additionally choose an edge of $[-1, 1]^4$, parallel to the $x_4$-axis, skew to the first three edges. Last we choose a vertex of $[-1, 1]^4$ with 4-th coordinate 1, at distance 1 from the last chosen edge. Then bend each edge to a circular arc as in the proof of Theorem 3.

# 6 Antipodal and strictly antipodal sets of segments

**Proof of Theorem 4.1:** We deal with the antipodal and strictly antipodal cases parallelly.

Observe that antipodality implies that the lines spanned by the segments are not intersecting or coincident, i.e., they are either parallel and distinct or skew.

First suppose that the lines spanned by the segments are pairwise skew. Then replace the segments by some segments in their relative interiors. Thus antipodality or strict antipodality of the original segments implies strict antipodality of the new segments. Then by [21, Proposition 2], i.e., by (11) of this paper, $n \leq 3$. This settles the antipodal case for pairwise skew segments. There remains the strictly antipodal case for pairwise skew segments.

We follow the proof of [21, Proposition 2]. We may suppose that we have $n = 3$ segments $s_1, s_2, s_3$. Then $[cl\, conv\, (s_1 \cup s_2)] \cap [cl\, conv\, (s_2 \cup s_3)] \cap [cl\, conv\, (s_3 \cup s_1)]$ is a parallelepiped, and $s_1, s_2, s_3$ lie on three mutually skew edge lines of it. By strict antipodality, $s_1, s_2, s_3$ lie actually in the relative interiors of the mutually skew edges of this parallelepiped.

Second suppose that all segments are parallel, but no two of them are collinear. Then, projecting them along their direction to a plane, we obtain an antipodal or strictly antipodal set of points in this plane, respectively. This implies $n \leq 4$, or $n \leq 3$, respectively, with equality only if the projections form the vertices of a parallelogram, or of a triangle, respectively. This proves the statement of the theorem for parallel segments.

Third suppose that, e.g., $s_1$ and $s_2$ are parallel, but not collinear, and not all other $s_i$'s are parallel to them. Then, for $i \neq 1, 2$, $\text{aff}\, s_i$ is either
(1) parallel to but distinct from $\text{aff}\, s_1$, $\text{aff}\, s_2$, and then by antipodality $\text{aff}\, s_1$, $\text{aff}\, s_2$, $\text{aff}\, s_3$ are not collinear, or
(2) $\text{aff}\, s_i$ is skew to $\text{aff}\, s_1$, $\text{aff}\, s_2$.
In case (2) $\text{aff}\, s_i$ cannot intersect the plane $\text{aff}\, (s_2 \cup s_3)$, either in $\text{conv}\, ((\text{aff}\, s_1) \cup (\text{aff}\, s_2))$, or outside of it, by antipodality.

In case(1), if $n = 3$, we have the statement of the theorem. In case (2) we do not have strict antipodality, only antipodality. Therefore we may assume $n \geq 4$. Both in case (1) and (2), $\text{aff}\, s_i$ is parallel to $\text{aff}\, (s_1 \cup s_2)$.

Let us suppose that among the $s_i$'s, for $i \neq 1, 2$, there is one parallel to $s_1, s_2$, say, $s_3$, and also one skew to $s_1, s_2$, say $s_4$. Then $\text{aff}\, s_1$, $\text{aff}\, s_2$, $\text{aff}\, s_3$ are not collinear, by antipodality. Then $\text{aff}\, s_4$ is parallel both to $\text{aff}\, (s_1 \cup s_2)$ and $\text{aff}\, (s_2 \cup s_3)$. Hence $s_4$ is parallel to $s_1, s_2$, a contradiction.

Therefore, the $s_i$'s, for $i \neq 1, 2$, are either all parallel to $s_1, s_2$, or all skew to $s_1, s_2$. The first case was settled above. Therefore we suppose that all $s_i$'s, for $i \neq 1, 2$, are skew to $s_1, s_2$.

Suppose that $\text{aff}\, (s_1 \cup s_2)$ is horizontal, and suppose that $\text{aff}\, s_3$, $\text{aff}\, s_4$ are skew to $\text{aff}\, s_1$, $\text{aff}\, s_2$. Then $\text{aff}\, s_3$, $\text{aff}\, s_4$ cannot be in different heights above $\text{aff}\, (s_1 \cup$ $s_2$), by antipodality. Hence aff $(s_3 \cup s_4)$ is parallel to and distinct from aff $(s_1 \cup s_2)$, and $s_3, s_4$ are parallel. If there were a third segment $s_3$ skew to $s_1, s_2$, then $s_3, s_4, s_5$ would be coplanar and parallel, contradicting antipodality. Since $s_1, s_2$, as well as $s_3, s_4$, are antipodal, we have the statement of the theorem, second case. □

**Proof of Corollary 4.2:** For an antipodal, or strictly antipodal set of lines in $\mathbb{R}^3$, any of their subsegments form an antipodal, or strictly antipodal set of segments in $\mathbb{R}^3$. Hence the inequalities in Theorem 4.1 imply the inequalities in Corollary 4.2.

Both for the antipodal, and strictly antipodal case, the first cases of equality in Theorem 4.1 imply the cases of equality in Corollary 4.2. There remains to show that, for the second cases of equality in Theorem 4.1, the lines spanned by the segments are not antipodal, or strictly antipodal.

The lines spanned by two parallel edges of a face of a parallelepiped lie on different boundary planes of parallel slabs which can be even rotated about the lines spanned by these edges. However, such a slab never contains the lines spanned by the two parallel edges of the opposite face, which are not parallel to the first mentioned edges.

From among the lines spanned by three mutually skew edges of a parallelepiped let us choose two ones — say, they are horizontal. They lie on the different boundary planes of a single parallel slab, namely the horizontal one. However this slab does not contain any non-horizontal line. □

**Acknowledgements:** The authors are indebted to Z. Füredi for his help with an earlier variant of this paper and calling the authors' attention to the paper [10]; to I. Bárány and P. McMullen for their help with an earlier variant of this paper; and to K. J. Swanepoel for his useful remarks on an earlier variant of the paper and calling the authors' attention to the paper [7].